\documentclass[10pt]{amsart}

\usepackage{amssymb}
\usepackage{diagrams}
\newtheorem{proposition}{Proposition}[section]
\newtheorem{lemma}[proposition]{Lemma}
\newtheorem{corollary}[proposition]{Corollary}
\newtheorem{theorem}[proposition]{Theorem}

\theoremstyle{definition}
\newtheorem{definition}[proposition]{Definition}
\newtheorem{example}[proposition]{Example}
\newtheorem{examples}[proposition]{Examples}

\theoremstyle{remark}
\newtheorem{remark}[proposition]{Remark}
\newtheorem{remarks}[proposition]{Remarks}

\newcommand{\thlabel}[1]{\label{th:#1}}
\newcommand{\thref}[1]{Theorem~\ref{th:#1}}
\newcommand{\selabel}[1]{\label{se:#1}}
\newcommand{\seref}[1]{Section~\ref{se:#1}}
\newcommand{\lelabel}[1]{\label{le:#1}}
\newcommand{\leref}[1]{Lemma~\ref{le:#1}}
\newcommand{\prlabel}[1]{\label{pr:#1}}
\newcommand{\prref}[1]{Proposition~\ref{pr:#1}}

\newcommand{\relabel}[1]{\label{re:#1}}

\newcommand{\exlabel}[1]{\label{ex:#1}}
\newcommand{\exref}[1]{Example~\ref{ex:#1}}
\newcommand{\exslabel}[1]{\label{ex:#1}}
\newcommand{\exsref}[1]{Example~\ref{ex:#1}}
\newcommand{\delabel}[1]{\label{de:#1}}
\newcommand{\deref}[1]{Definition~\ref{de:#1}}
\newcommand{\eqlabel}[1]{\label{eq:#1}}
\newcommand{\equref}[1]{(\ref{eq:#1})}

\def\NN{{\mathbb N}}

\def\equal#1{\smash{\mathop{=}\limits^{#1}}}

\def\<{\langle}
\def\>{\rangle}

\def\va{\varepsilon}
\def\v{\varphi}

\def\tr{\triangleright}
\def\rh{\rightharpoonup}
\def\lh{\leftharpoonup}

\def\ra{\rightarrow}
\def\a{\alpha}
\def\b{\beta}

\def\l{\lambda}

\def\cd{\cdot}
\def\d{\delta}
\def\O{\Omega}

\def\ov{\overline}
\def\un{\underline}
\def\mf{\mathfrak}

\def\le{\langle}
\def\ri{\rangle}
\newcommand{\tpla}{\mbox{$\tilde {p}^1$}}
\newcommand{\tplb}{\mbox{$\tilde {p}^2$}}

\newcommand{\tqla}{\mbox{$\tilde {q}^1$}}
\newcommand{\tqlb}{\mbox{$\tilde {q}^2$}}

\newcommand{\smi}{\mbox{$S^{-1}$}}

\def\rawo\lonra{\longrightarrow}

\def\ot{\otimes}

\def\cal{\mathcal}

\def\equal#1{\smash{\mathop{=}\limits^{#1}}}

\begin{document}
\title[Involutory quasi-Hopf algebras]
{Involutory quasi-Hopf algebras}
\author{D. Bulacu}
\address{Faculty of Mathematics and Informatics, University of Bucharest,
Str. Academiei 14, RO-70109, Bucharest 1, Romania}
\email{daniel.bulacu@fmi.unibuc.ro}
\author{S. Caenepeel}
\address{Faculty of Engineering, 
Vrije Universiteit Brussel, B-1050 Brussels, Belgium}
\email{scaenepe@vub.ac.be}
\urladdr{http://homepages.vub.ac.be/\~{}scaenepe/}
\author{B. Torrecillas}
\address{Department of Algebra and Analysis\\
Universidad de Almer\'{\i}a\\
E-04071 Almer\'{\i}a, Spain}
\email{btorreci@ual.es}
\thanks{Research supported by the EU research program LIEGRITS, 
RTN 2003, 505078, and the bilateral project BWS04/04 ``New techniques in 
Hopf algebras and graded ring theory" of the Flemish and Romanian 
Governments, and the grant 434/1.10.2007 of CNCSIS (project ID 1005).}
\thanks{
This paper was written while the first author was
visiting the Universidad de Almer\'{\i}a
and the Vrije Universiteit Brussel. 
He would like to thank both universities for their warm hospitality.}
\dedicatory{To Fred Van Oystaeyen on the occasion of his 60th birthday.}

\begin{abstract}
We introduce and investigate the basic properties of 
an involutory (dual) quasi-Hopf algebra. We also study the representations 
of an involutory quasi-Hopf algebra and prove that an involutory 
dual quasi-Hopf algebra with non-zero integral is cosemisimple.  
\end{abstract}
\maketitle
\section*{Introduction}\label{sec0}
One of the aims in \cite{bt} was to find a plausible definition 
of an involutory quasi-Hopf algebra. Since the definition 
of a quasi-Hopf algebra $H$ is given in such a way that the category 
of its finite dimensional representations ${}_H{\cal M}^{\rm fd}$ has a 
rigid monoidal structure it seems to be natural to relate the involutory 
notion to a certain property of this category.

On one hand, in the 
classical case of a Hopf algebra, a categorical interpretation for 
the involutory notion was given by Majid in \cite{m2}. Namely, he has shown 
that for a finite dimensional Hopf algebra $H$, the trace of the square of the antipode 
$S$ of $H$, ${\rm Tr}(S^2)$, arises in a very natural way as the representation-theoretic
rank of the Schr${\rm{\ddot{o}}}$dinger representation of $H$, $\un{\rm dim}(H)$, or as the
representation-theoretic rank of the canonical representation of the quantum double,
$\un{\rm dim}(D(H))$. Relating this to results of Larson and Radford \cite{lr1, lr2} 
and Etingof and Gelaki \cite{eg2}, we obtain that the above rank coincide to the classical dimension 
of $H$ if $H$ is involutory; otherwise it is zero.

This is why, for a finite dimensional quasi-Hopf algebra $H$, we 
computed in \cite{bt} the representation-theoretic rank of a finite dimensional 
quasi-Hopf algebra $H$, and of its associated quantum double $D(H)$, 
within the category of 
finite dimensional left modules over $D(H)$, ${}_{D(H)}{\cal M}^{\rm fd}$. 
More precisely, for a quasi-Hopf algebra $H$, we have 
\[
\un{\rm dim}(H)=\un{\rm dim}(D(H))={\rm Tr}\left(h\mapsto S^{-2}(S(\b)\a h\b S(\a))\right).
\]
Therefore, we call a quasi-Hopf algebra $H$ involutory if $S^2(h)=S(\b)\a h\b S(\a)$, 
for all $h\in H$.

On the other hand, due to some recent results of Etingof, Nikshych 
and Ostrik \cite{eno}, it seems that the notion of involutory quasi-Hopf algebra 
should be given in such a way that, at least over an algebraic closed field 
of characteristic zero, its category of finite dimensional 
representations is a fusion category and, moreover, has that unique pivotal structure 
with respect to which the categorical 
dimensions of simple objects coincide with their usual dimensions. 
In other words, an involutory quasi-Hopf algebra, say $H$, should be semisimple as 
an algebra and such that 
the identity functor and the second duality functor 
$(-)^{**}$ of the category ${}_H{\cal M}^{\rm fd}$ 
are tensor isomorphic via a tensor functor, say $j$. In addition, we should have 
\[
{\rm dim}_k(V)={\rm ev}_{V^*}\circ (j_V\ot {\rm id}_{V^*})\circ {\rm coev}_V,
\]
for any simple finite dimensional left $H$-module $V$, where ${\rm ev}_{V^*}$ and ${\rm coev}_V$ 
are the evaluation map of $V^*$ and the coevaluation map of $V$, respectively. Note that 
the composition on the right hand side of the above equality was defined in \cite{eno} 
as being the categorical dimension of $V$.

We have to stress the fact that our definition for 
an involutory quasi-Hopf algebra agrees with this point of 
view. More precisely, over a field of characteristic zero (or, more generally, 
if ${\rm dim}(H)\not=0$ in $k$) 
any involutory quasi-Hopf algebra is semisimple because of 
the trace formula for quasi-Hopf algebras proved in \cite{bt}. Also, by 
\leref{2.2} below for an involutory quasi-Hopf algebra 
the square of the antipode in an inner automorphism 
defined by $S(\b)\a$. Now, using some results from 
\cite{mn, sch} we will prove that the invertible element ${\mf g}$ which defines $j$ is 
exactly $\b S(\a)$, the inverse of $S(\b)\a$. From here we conclude that 
the family of left $H$-module isomorphisms   
\[
j_V:\ V\in v\mapsto (v^*\mapsto v^*(S(\b)\a \cd v))\in V^{**}
\] 
endows ${}_H{\cal M}^{\rm fd}$ with the unique pivotal structure 
with respect to which the categorical dimension 
of any simple object $V$ of ${}_H{\cal M}^{\rm fd}$ coincides with its usual dimension. 
All the details are presented in \seref{pivotal}.
  
Then the aim of this paper is to see that with our definition of an involutory (dual) 
quasi-Hopf algebra some results of Lorenz \cite{lz} can be generalized 
to quasi-Hopf algebras, 
and that the main result of D\u asc\u alescu, N\u ast\u asescu and Torrecillas 
in \cite{dnt} has a counterpart for involutory dual quasi-Hopf algebras.

The paper is organized as follows. 
In \seref{2} we introduce the involutory notion and then we study their 
basic properties. Starting with $H(2)$, the (involutory) quasi-Hopf algebra 
presented in \cite{eg}, using different constructions as transmutation, bosonozation or 
quantum double, we will be able to construct three involutory quasi-Hopf algebras 
of dimension $4$. It comes out that all of them are $k[C_2\times C_2]$, 
the group Hopf algebra associated to the Klein group, viewed as quasi-Hopf algebras via some 
non-equivalent $3$-cocycles. 
Moreover, all of them are not coming from twisting $k[C_2\ot C_2]$ 
by a gauge transformation. Also, we will prove in \seref{pivotal} that, under a ``natural" condition, 
the involutory property survives when we pass from $H$ to 
its quantum double. Note that this ``natural" condition is satisfied if we work over 
an algebraic closed field of characteristic zero.

In \seref{3} we study the representations of an involutory quasi-Hopf algebra. 
Namely, we show that if $H$ is an involutory semisimple quasi-Hopf algebra 
over a field $k$ then the 
characteristic of $k$ does not divide the dimension of any finite dimensional 
absolutely simple $H$-module. We also prove that if $H$ is a non-semisimple 
involutory quasi-Hopf algebra then the  characteristic 
of $k$ divides the dimension of any finite dimensional projective $H$-module. Both 
results generalize well-known results in Hopf algebra theory 
due to Larson and Lorenz, see \cite{l, lz}.

The first goal of \seref{4} is to introduce and study the dual concept: 
involutory dual quasi-Hopf algebras. Then due to some recent results 
proved in \cite{bc1} we will be able to show that any involutory co-Frobenius dual 
quasi-Hopf algebra over a field of characteristic zero is cosemisimple. 
This result was proved for Hopf algebras in \cite{dnt} and it can be viewed as an 
extension of the results of Larson and Radford to the infinite dimensional case. 
For Hopf algebras it is well-known that the converse property is true only in 
the finite dimensional case.

\section{Preliminaries}\selabel{1}
We work over a commutative field $k$. All algebras, linear
spaces etc. will be over $k$; unadorned $\ot $ means $\ot_k$.
Following Drinfeld \cite{d1}, a quasi-bialgebra is
a four-tuple $(H, \Delta , \va , \Phi )$ where $H$ is
an associative algebra with unit,
$\Phi$ is an invertible element in $H\ot H\ot H$, and
$\Delta :\ H\ra H\ot H$ and $\va :\ H\ra k$ are algebra
homomorphisms satisfying the identities
\begin{eqnarray}
&&(id \ot \Delta )(\Delta (h))=
\Phi (\Delta \ot id)(\Delta (h))\Phi ^{-1},\label{q1}\\
&&(id \ot \va )(\Delta (h))=h,
\mbox{${\;\;\;}$}
(\va \ot id)(\Delta (h))=h,\label{q2}
\end{eqnarray}
for all $h\in H$, and
$\Phi$ has to be a $3$-cocycle, in the sense that
\begin{eqnarray}
&&(1\ot \Phi)(id\ot \Delta \ot id)
(\Phi)(\Phi \ot 1)=(id\ot id \ot \Delta )(\Phi )(\Delta \ot id \ot
id)(\Phi ),\label{q3}\\
&&(id \ot \va \ot id)(\Phi )=1\ot 1.\label{q4}
\end{eqnarray}
The map $\Delta$ is called the coproduct or the
comultiplication, $\va $ the counit and $\Phi $ the
reassociator. As for Hopf algebras we denote $\Delta (h)=h_1\ot h_2$,
but since $\Delta $ is only quasi-coassociative we adopt the
further convention (summation understood):
$$
(\Delta \ot id)(\Delta (h))=h_{(1, 1)}\ot h_{(1, 2)}\ot h_2,
\mbox{${\;\;\;}$}
(id\ot \Delta )(\Delta (h))=h_1\ot h_{(2, 1)}\ot h_{(2,2)},
$$
for all $h\in H$. We will
denote the tensor components of $\Phi$
by capital letters, and the ones of
$\Phi^{-1}$ by small letters,
namely
\begin{eqnarray*}
&&\Phi=X^1\ot X^2\ot X^3=T^1\ot T^2\ot T^3=
V^1\ot V^2\ot V^3=\cdots\\
&&\Phi^{-1}=x^1\ot x^2\ot x^3=t^1\ot t^2\ot t^3=
v^1\ot v^2\ot v^3=\cdots
\end{eqnarray*}
$H$ is called a quasi-Hopf
algebra if, moreover, there exists an
anti-morphism $S$ of the algebra
$H$ and elements $\a , \b \in
H$ such that, for all $h\in H$, we
have:
\begin{eqnarray}
&&S(h_1)\a h_2=\va(h)\a
\mbox{${\;\;\;}$ and ${\;\;\;}$}
h_1\b S(h_2)=\va (h)\b,\label{q5}\\ 
&&X^1\b S(X^2)\a X^3=1
\mbox{${\;\;\;}$ and${\;\;\;}$}
S(x^1)\a x^2\b S(x^3)=1.\label{q6}
\end{eqnarray}
Our definition of a quasi-Hopf algebra is different from the
one given by Drinfeld \cite{d1} in the sense that we do not
require the antipode to be bijective. Nevertheless, in the finite
dimensional or quasi-triangular case this condition can be deleted because
it follows from the other axioms, see \cite{bc1} and \cite{bn3}.

The definition of a quasi-Hopf algebra is ``twist coinvariant" in the following sense.
An invertible element $F\in H\ot H$ is called a gauge 
transformation or twist if
$(\va \ot id)(F)=(id\ot \va)(F)=1$. 
If $H$ is a quasi-Hopf algebra and $F=F^1\ot F^2\in H\ot H$ 
is a gauge transformation with inverse $F^{-1}=G^1\ot G^2$, then we can define a 
new quasi-Hopf algebra $H_F$ by keeping 
the multiplication, unit, counit 
and antipode of $H$ and 
replacing the comultiplication, reassociator and the 
elements $\alpha$ and $\beta$ by 
\begin{eqnarray}
&&\Delta _F(h)=F\Delta (h)F^{-1},\label{g1}\\
&&\Phi_F=(1\ot F)(id \ot \Delta )(F) \Phi (\Delta \ot
id)(F^{-1})(F^{-1}\ot 1),\label{g2}\\
&&\a_F=S(G^1)\a G^2,
\mbox{${\;\;\;}$}%
\b_F=F^1\b S(F^2).\label{g3}
\end{eqnarray}
For a quasi-Hopf algebra the antipode is
determined uniquely up to a transformation 
$\a \mapsto \a _{\mathbb{U}}:=\mathbb{U}\a $,
$\b \mapsto \b _{\mathbb{U}}:=\b \mathbb{U}^{-1}$,
$S(h)\mapsto S_{\mathbb{U}}(h):=\mathbb{U}S(h)\mathbb{U}^{-1}$,
where $\mathbb{U}\in H$ is invertible. In this case we will denote
by $H^{\mathbb{U}}$ the new quasi-Hopf algebra $(H, \Delta , \va, \Phi ,
S_{\mathbb{U}}, \a _{\mathbb{U}}, \b _{\mathbb{U}})$. \\
${\;\;\;}$
If $H=(H, \Delta , \va , \Phi , S, \a , \b )$ is 
a quasi-bialgebra or a quasi-Hopf algebra then 
$H^{\rm op}$, $H^{\rm cop}$ and $H^{\rm op, cop}$ are also quasi-bialgebras 
(respectively quasi-Hopf algebras), 
where "op" means opposite 
multiplication and "cop" means opposite comultiplication. The 
structures are obtained by putting $\Phi _{\rm op}=\Phi ^{-1}$, 
$\Phi _{\rm cop}=(\Phi ^{-1})^{321}$, $\Phi _{\rm op, cop}=\Phi ^{321}$, 
$S_{\rm op}=S_{\rm cop}=(S_{\rm op, cop})^{-1}=S^{-1}$, $\a _{\rm op}=\smi (\b )$, 
$\b _{\rm op}=\smi (\a )$, $\a _{\rm cop}=\smi (\a )$, $\b _{\rm cop}=\smi (\b )$, 
$\a _{\rm op, cop}=\b $ and $\b _{\rm op, cop}=\a $.

The axioms for a quasi-Hopf algebra imply that
$\va \circ S=\va $ and $\va (\a )\va (\b )=1$,
so, by rescaling $\a $ and $\b $, we may assume without loss of generality
that $\va (\a )=\va (\b )=1$. The identities
(\ref{q2}), (\ref{q3}) and (\ref{q4}) also imply that
\begin{equation}\label{q7}
(\va \ot id\ot id)(\Phi )=
(id \ot id\ot \va )(\Phi )=1\ot 1.
\end{equation}
It is well-known that the antipode of a Hopf algebra is an
anti-coalgebra morphism. For a quasi-Hopf algebra, we have
the following statement: there exists a gauge transformation
$f\in H\ot H$ such that
\begin{equation} \label{ca}
f\Delta (S(h))f^{-1}= (S\ot S)(\Delta ^{\rm cop}(h))
\mbox{,${\;\;\;}$for all $h\in H$,}
\end{equation}
where $\Delta ^{\rm cop}(h)=h_2\ot h_1$. The element $f$ can be computed
explicitly. First set
\begin{equation}
A^1\ot A^2\ot A^3\ot A^4=
(\Phi \ot 1) (\Delta \ot id\ot id)(\Phi ^{-1}),
\end{equation}
\begin{equation} 
B^1\ot B^2\ot B^3\ot B^4=
(\Delta \ot id\ot id)(\Phi )(\Phi ^{-1}\ot 1)
\end{equation}
and then define $\gamma, \delta\in H\ot H$ by
\begin{equation} \label{gd}
\gamma =S(A^2)\a A^3\ot S(A^1)\a A^4~~{\rm and}~~
\delta
=B^1\b S(B^4)\ot B^2\b S(B^3).
\end{equation}
With this notation $f$ and $f^{-1}$ are given by the formulas
\begin{eqnarray}
f&=&(S\ot S)(\Delta ^{\rm op}(x^1)) \gamma \Delta (x^2\b
S(x^3)),\label{f}\\ 
f^{-1}&=&\Delta (S(x^1)\a x^2) \delta
(S\ot S)(\Delta^{\rm cop}(x^3)).\label{g}
\end{eqnarray}
Moreover, $f$ satisfies the following relations:
\begin{equation} \label{gdf}
f\Delta (\a )=\gamma ,
\mbox{${\;\;\;}$}
\Delta (\b )f^{-1}=\delta ,
\end{equation}
\begin{equation} \label{pf}
(1\ot f)(id \ot \Delta )(f) \Phi (\Delta \ot id)
(f^{-1})(f^{-1}\ot 1)
=(S\ot S\ot S)(X^3\ot X^2\ot X^1).
\end{equation}
${\;\;\;}$
In a Hopf algebra $H$, we obviously have the identity
\[
h_1\ot h_2S(h_3)=h\ot 1,~{\rm for~all~}h\in H.
\]
We will need the generalization of this formula to
quasi-Hopf algebras. Following \cite{hn1, hn2}, we define
\begin{eqnarray}
&&p_R=p^1\ot p^2=x^1\ot x^2\b S(x^3),~~
q_R=q^1\ot q^2=X^1\ot S^{-1}(\a X^3)X^2,\label{qr}\\
&&p_L=\tpla \ot \tplb=X^2\smi (X^1\b )\ot X^3,~~
q_L=\tqla \ot \tqlb=S(x^1)\a x^2\ot x^3.\label{ql}
\end{eqnarray}
For all $h\in H$, we then have
\begin{eqnarray}
\Delta (h_1)p_R(1\ot S(h_2))&=&p_R(h\ot 1),\label{qr1}\\
(1\ot S^{-1}(h_2))q_R\Delta (h_1)&=&(h\ot 1)q_R,\label{qr1a}\\
\Delta (h_2)p_L(\smi (h_1)\ot 1)&=&p_L(1\ot h),\label{ql1}\\
(S(h_1)\ot 1)q_L\Delta (h_2)&=&(1\ot h)q_L.\label{ql1a}
\end{eqnarray}
Furthermore, the following relations hold
\begin{eqnarray}
&&(1\ot S^{-1}(p^2))q_R\Delta (p^1)=1\ot 1,\label{pqra}\\
&&\Delta (q^1)p_R(1\ot S(q^2))=1\ot 1,\label{pqr}\\
&&(S(\tpla)\ot 1)q_L\Delta (\tplb)=1\ot 1,\label{pql}\\
&&\Delta (\tqlb)p_L(\smi (\tqla)\ot 1)=1\ot 1.\label{pqla}
\end{eqnarray}
Finally, for further use we need the notion of quasi-triangular quasi-Hopf 
algebra. Recall that a quasi-Hopf algebra $H$ is
quasi-triangular ($QT$ for short) if there exists an element
$R=R^1\ot R^2=r^1\ot r^2\in H\ot H$ such that
\begin{eqnarray}
(\Delta \ot id)(R)&=&X^2R^1x^1Y^1\ot X^3x^3r^1Y^2\ot
X^1R^2x^2r^2Y^3\label{qt1}\\ 
(id \ot \Delta )(R)&=&x^3R^1X^2r^1y^1\ot x^1X^1r^2y^2\ot 
x^2R^2X^3y^3\label{qt2}\\
\Delta ^{\rm cop}(h)R&=&R\Delta (h),~{\rm for~all~}h\in
H\label{qt3}\\ 
(\va \ot id)(R)&=&(id\ot \va)(R)=1.\label{qt4}
\end{eqnarray}
To any finite dimensional quasi-Hopf algebra $H$ we can associate 
a $QT$ quasi-Hopf algebra $D(H)$, the quantum double of $H$. 
From \cite{hn1, hn2, bc2}, we recall the definition of the quantum double
$D(H)$. Let $\{e_i\}_{i=\ov {1, n}}$ be a basis of $H$, and 
$\{e^i\}_{i=\ov{1, n}}$ the corresponding dual basis of $H^*$. 
We can easily see that $H^*$, 
the linear dual of $H$, is not a quasi-Hopf algebra.
But $H^*$ has a dual structure coming from the initial structure of $H$.
So $H^*$ is a coassociative coalgebra, with comultiplication
\[
\widehat {\Delta }(\v )=\v _1\ot \v _2=
\sum \limits _{i, j=1}^n\v (e_ie_j)e^i\ot e^j,
\]
or, equivalently, $\widehat {\Delta }(\v )=\v _1\ot \v_2 
\Leftrightarrow \v (hh')=\v _1(h)\v _2(h')$, for all $h, h'\in H$.

$H^*$ is also an $H$-bimodule, by
\[
\<h\rh \v , h'\>=\le \v ,h'h\ri,
\mbox{${\;\;\;}$}
\<\v \lh h, h'\>=\le\v ,hh'\ri.
\]
The convolution is a multiplication on $H^*$; it is not
associative, but only quasi-associative:
\[
[\v \psi]\xi=(X^1\rh \v \lh x^1)[(X^2\rh \psi \lh x^2)
(X^3\rh \xi \lh x^3)],~~\forall ~~\v , \psi , \xi \in H^*.
\]
We also introduce $\ov {S}:\ H^*\ra H^*$ as the coalgebra
antimorphism dual to $S$, this means
$\<\ov {S}(\v), h\>=$ $\<\v, S(h)\>$, for all $\v \in H^*$ and $h\in H$.

Now consider $\Omega \in H^{\ot 5}$ given by
\begin{eqnarray}
&&\hspace*{-5mm} \Omega =\Omega^1\ot \Omega^2\ot \Omega^3\ot
\Omega^4\ot \Omega^5 \nonumber \\
&&\hspace*{5mm}
=X^1_{(1, 1)}y^1x^1\ot X^1_{(1, 2)}y^2x^2_1\ot
X^1_2y^3x^2_2\ot \smi (f^1X^2x^3)\ot \smi (f^2X^3), \label{O}
\end{eqnarray}
where $f\in H\ot H$ is the element defined in (\ref{f}). We define
the quantum double $D(H)=H^*\bowtie H$ as follows: as a $k$-linear
space, $D(H)$ equals $H^*\ot H$, and the multiplication is given
by
\begin{eqnarray}
&&\hspace*{-1cm}(\v \bowtie h)(\psi \bowtie h')\nonumber\\
&=&[(\Omega ^1\rh \v \lh \Omega ^5)(\Omega
^2h_{(1, 1)}\rh \psi \lh \smi (h_2)\Omega ^4)]\bowtie \Omega
^3h_{(1, 2)}h'.\label{mdd}
\end{eqnarray}

From \cite{hn1,hn2} we know that $D(H)$ is an
associative algebra with unit $\va \bowtie 1$, and $H$ is a unital
subalgebra via the morphism $i_D:\ H\ra D(H)$, $i_D(h)=\va \bowtie
h$. Moreover, $D(H)$ is a quasi-triangular quasi-Hopf algebra
with the following structure:
\begin{eqnarray}
&&\hspace*{-5mm}\Delta _D(\v \bowtie h)=
(\va \bowtie X^1Y^1)
(p^1_1x^1\rh \v _2\lh Y^2\smi (p^2)\bowtie p^1_2x^2h_1)\nonumber\\
&&\hspace*{10mm}\ot (X^2_1\rh \v _1\lh \smi (X^3)\bowtie X^2_2Y^3x^3h_2)
\label{cddf}\\
&&\hspace*{-5mm}
\va _D(\v \bowtie h)=\va (h)\v (\smi (\a ))\label{coddf}\\
&&\hspace*{-5mm}
\Phi _D=(i_D\ot i_D\ot i_D)(\Phi)\label{pdd}\\
&&\hspace*{-5mm}
S_D(\v \bowtie h)=(\va \bowtie S(h)f^1)(p^1_1U^1\rh \ov
{S}^{-1}(\v )\lh f^2\smi (p^2)\bowtie p^1_2U^2)\label{anddf}\\
&&\hspace*{-5mm}
\a _D=\va \bowtie \a ,~~\b _D=\va \bowtie \b \label{abdd}\\
&&\hspace*{-5mm}
R_D=\sum \limits _{i=1}^n
(\va \bowtie \smi (p^2)e_ip^1_1)\ot (e^i\bowtie p^1_2).\label{rdd}
\end{eqnarray}
Here $p_R=p^1\ot p^2$ and $f=f^1\ot f^2$ are the elements defined by
(\ref{qr}) and (\ref{f}), respectively, and $U=U^1\ot U^2\in H\ot H$
is the following element
\begin{equation}\label{U}
U=U^1\ot U^2=g^1S(q^2)\ot g^2S(q^1),
\end{equation}
where $f^{-1}=g^1\ot g^2$ and $q_R=q^1\ot q^2$ are the elements
defined by (\ref{g}) and (\ref{qr}), respectively.
\section{Involutory quasi-Hopf algebras}\selabel{2}
\setcounter{equation}{0}
The purpose of this Section is to introduce and study involutory quasi-Hopf 
algebras. As we have already explained in the Introduction there is a categorical 
interpretation for this definition, see \cite{bt}. 
\begin{definition}\delabel{2.1}
A quasi-Hopf algebra is called involutory if the following formula holds for
all $h\in H$:
\begin{equation}\label{inv}
S^2(h)=S(\b)\a h\b S(\a).
\end{equation}
\end{definition}

Examples of involutory quasi-Hopf algebras will be presented at the end of this 
Section. 

Next we prove that for an involutory quasi-Hopf algebra the square of 
the antipode is inner via an element depending of $\a$ and $\b$.

\begin{lemma}\lelabel{2.2}
Let $H$ be an involutory quasi-Hopf algebra. Then $S(\b)\a$ is 
an invertible element and $(S(\b)\a)^{-1}=\b S(\a)$. In particular, $S^2$ is inner 
and therefore $S$ is bijective. Moreover, $\a$ and $\b$ 
are invertible elements and 
\begin{eqnarray}
&&\a ^{-1}=\smi (\a \b)\b =\b S(\b \a),\label{inva}\\
&&\b ^{-1}=S(\b \a)\a =\a \smi (\a \b).\label{invb}  
\end{eqnarray}
\end{lemma}

\begin{proof}
For simplicity denote $u =S(\b)\a$ and $v  =\b S(\a)$. Then 
$S^2(h)=u hv  $, for all $h\in H$. Since $S^2$ is an algebra map, we have
that $1=S^2(1)=u v  $. This implies that $S^2(u )=u u v  =u $
and $S^2(v  )=u v v =v$. Then we find that 
$v u =S^2(v  )S^2(u )=S^2(v  u )=u v u v=1$. 
This shows that $u^{-1}=v$. The relation $vu =1$ comes out as 
$\b S(\b \a)\a =1$. Thus $\b$ has a 
right inverse, namely $S(\b \a)\a $,
and $\alpha$ has a left inverse, namely $\b S(\b \a)$. Similarly, from $u v  =1$ 
we obtain that $S(\b)\a \b S(\a)=1$. Since $S$ is bijective this is equivalent to 
$\a \smi (\a \b)\b =1$, so $\b $ has also a left inverse, namely $\a \smi (\a \b)$,
and $\alpha$ has a right inverse, namely $\smi (\a \b)\b$. It 
follows now that $\alpha $ and $\beta$ are invertible and that there inverse are
given by (\ref{inva},\ref{invb}).
\end{proof}

\begin{remark}\label{re2.3}
Following the ideas in \cite{bt}, there are two seemingly different ways to introduce the
notion of involutory quasi-Hopf algebra. The first definition is obtained from the
formula of the categorical representation rank of $H$ and $D(H)$, the quantum double of $H$, and
requires that the map sending $h$ to $S^{-2}(S(\b)\a h\b S(\a))$ is the identity of $H$;
this is clearly equivalent to (\ref{inv}) in \deref{2.1}. The second definition involves $S^2$
and is obtained from the trace formula for quasi-Hopf algebras proved in \cite{bt}. It
requires that the map sending $h\in H$ to $\b S(\a)S^2(h)S(\b)\a$ is the identity of $H$.
It follows immediately from \leref{2.2} that these two definitions are equivalent.
Moreover, we can easily verify that $H$ is involutory if 
and only if $H^{\rm op}$ is involutory, if and only if $H^{\rm cop}$ 
is involutory, if and only if $H^{\rm op, cop}$ is involutory.
\end{remark}

For a Hopf algebra $H$ it is well-known that $S^2=id_H$ if and only if 
$S(h_2)h_1=\va(h)1$ for any $h\in H$, if and only if $h_2S(h_1)=\va (h)1$ 
for any $h\in H$ (see for instance \cite[Proposition 4.2.7]{dnr}). For 
quasi-Hopf algebras we have the following result. 

\begin{proposition}\prlabel{2.4}
Let $H$ be an involutory quasi-Hopf algebra. Then for all $h\in H$ the following 
relations hold:
\begin{equation}\label{sp}
S(h_2)\b ^{-1}h_1=\va (h)\b ^{-1}~~{\rm and}~~
h_2\a ^{-1}S(h_1)=\va (h)\a ^{-1}.
\end{equation}

\begin{proof}
As we have seen, if $\mathbb{U}\in H$ is invertible then we can define a new 
quasi-Hopf algebra $H^{\mathbb{U}}=(H, \Delta , \va, \Phi ,
S_{\mathbb{U}}, \a _{\mathbb{U}}, \b _{\mathbb{U}})$, where 
$\a _{\mathbb{U}}=\mathbb{U}\a $, $\b _{\mathbb{U}}=\b \mathbb{U}^{-1}$ and 
$S_{\mathbb{U}}(h)=\mathbb{U}S(h)\mathbb{U}^{-1}$.

Now, consider $\mathbb{U}=\a^{-1}$. We know from \leref{2.2} that $\mathbb{U}$ is 
invertible, so it makes sense to consider the quasi-Hopf algebra $H^{\mathbb{U}}$. 
In this particular case we have that 
$\a _{\mathbb{U}}=1$, $\b _{\mathbb{U}}=\b \a$ and 
\[
S_{\mathbb{U}}(h)=\a ^{-1}S(h)\a =\b S(\b \a)S(h)\a =\b S(\a)S(\b ^{-1}h\b )S(\b)\a=
\smi (\b ^{-1}h\b ).
\]  
Since $S_{\mathbb{U}}(h_1)\a _{\mathbb{U}}h_2=\va (h)\a _{\mathbb{U}}$ 
for all $h\in H$, we get that $\smi (\b ^{-1}h_1\b )h_2=\va (h)1$, 
and this is equivalent to 
$S(h_2)\b ^{-1}h_1\b =\va (h)1$, for all $h\in H$. It follows now that 
$S(h_2)\b ^{-1}h_1=\va (h)\b ^{-1}$ for all $h\in H$, as needed.

Similarly, using the fact that 
$h_1\b _{\mathbb{U}}S(h_2)=\va (h)\b _{\mathbb{U}}$ for all $h\in H$, 
one can prove that $h_2\a ^{-1}S(h_1)=\va (h)\a ^{-1}$ for all $h\in H$, 
the details are left to the reader. 
\end{proof}
\end{proposition}

Let us now present some examples of involutory quasi-Hopf algebras.

\begin{examples}\exslabel{exsinv}
1) To any finite dimensional cocommutative Hopf algebra $H$ 
and any normalized $3$-cocycle $\omega$ on $H$ we can associate a quasi-Hopf algebra, 
$D^{\omega}(H)$, see \cite{bp}. For $D^{\omega}(H)$ we have $\a=1$, $\b ^{-1}=S(\b)$ and 
$S^2(h)=\b^{-1}h\b$, for all $h\in D^{\omega}(H)$. Thus $D^{\omega}(H)$ is an 
involutory quasi-Hopf algebra.

2) Let $k$ be a field of characteristic different from $2$, 
$C_2$ the cyclic group of order two, and $g$ the generator of $C_2$. Following 
\cite{eg}, the two dimensional quasi-Hopf algebra $H(2)$ is the bialgebra 
$k[C_2]$, the group algebra associated to $C_2$, viewed as a quasi-Hopf algebra 
via the non-trivial reassociator $\Phi=1-2p_{-}\ot p_{-}\ot p_{-}$, 
where $p_{-} =\frac{1}{2}(1 - g)$. The antipode $S$ is the identity map and the 
distinguished elements $\a$ and $\b$ are $\a=g$ and $\b=1$, respectively. 
It is well-known that $H(2)$ is not twist equivalent to a Hopf algebra. 
It is easy to see that $H(2)$ is an involutory quasi-Hopf algebra.  
\end{examples}

We will now study further properties of $H(2)$. 
First we will show that if $k$ contains a primitive fourth 
root of unity then there are exactly two quasi-triangular 
structures on $H(2)$. In what follows, $p_{\pm}=\frac{1}{2}(1\pm g)$. One can easily check 
that $\{p_{\pm}\}$ is a basis for $H(2)$ consisting of orthogonal idempotents, and that 
$p_{+} + p_{-}=1$, $p_{+} - p_{-}=g$. 

\begin{proposition}\prlabel{2.6}
Suppose that $k$ is a field of characteristic different from $2$ containing
a primitive fourth root of unity $i$. Then there are exactly two different 
$R$-matrices for $H(2)$, namely, $R_{\pm}=1 -(1\pm i)p_{-}\ot p_{-}$.
\end{proposition}

\begin{proof}
Suppose that $R=a1\ot 1 + b1\ot g + cg\ot 1 + dg\ot g$ is an $R$-matrix for $H(2)$, 
where $a, b, c, d\in k$. By (\ref{qt4}) we have that 
$a + c=a + b=1$ and $b + d=c + d=0$, and therefore 
$b=c=-d$ and $a=1 - b$. Hence, $R$ should be on the form 
\[
R=(1 -b)1\ot 1 + b1\ot g + bg\ot 1 - bg\ot g=1 - b(1 - g)\ot (1 - g)
=1 - \omega p_{-}\ot p_{-},
\]
where we denoted $4b=\omega$. Now, we can easily see that
\begin{equation}\label{phih2}
\Phi ^{-1}=\Phi =1 - 2p_{-}\ot p_{-}\ot p_{-},
\end{equation} 
and since $X^2\ot X^3\ot X^1=\Phi$, the above relation implies 
\[
X^2R^1x^1\ot X^3x^3\ot X^1R^2x^2=1 - \omega p_{-}\ot 1\ot p_{-},
\] 
and therefore, after some computations, we get 
\begin{eqnarray*}
&&\hspace*{-1cm}
X^2R^1x^1Y^1\ot X^3x^3r^1Y^2\ot X^1R^2x^2r^2Y^3\\
&=&
1 -\omega p_{-}\ot p_{+}\ot p_{-} - \omega p_{+}\ot p_{-}\ot p_{-} 
- (2 - 2\omega + \omega^2)p_{-}\ot p_{-}\ot p_{-}.
\end{eqnarray*}
On the other hand, we have $\Delta(p_{-})=p_{-}\ot p_{+} + p_{+}\ot p_{-}$, so 
\[
(\Delta \ot id)(R)=1 -\omega p_{-}\ot p_{+}\ot p_{-} - \omega p_{+}\ot p_{-}\ot p_{-}.
\]
We conclude that (\ref{qt1}) holds if and only if $2 - 2\omega + \omega^2=0$, and 
this is equivalent to $\omega =1\pm i$.

Using (\ref{phih2}) for $\Phi^{-1}$, we obtain in a similar way that
\[
x^3R^1X^2\ot x^1X^1\ot x^2R^2X^3=1 - \omega p_{-}\ot 1\ot p_{-}.
\]
Using this formula, it can be proved that 
\begin{eqnarray*}
&&\hspace*{-1cm}
x^3R^1X^2r^1y^1\ot x^1X^1r^2y^2\ot x^2R^2X^3y^3\\
&&
=1 - \omega p_{-}\ot p_{-}\ot p_{+} - \omega p_{-}\ot p_{+}\ot p_{-} 
- (2 - 2\omega + \omega^2)p_{-}\ot p_{-}\ot p_{-}.
\end{eqnarray*}
It is easy to see that  
\[
(id\ot \Delta)(R)=1 - \omega p_{-}\ot p_{-}\ot p_{+} - \omega p_{-}\ot p_{+}\ot p_{-},
\]
so the relation in (\ref{qt2}) holds if and only if $2 - 2\omega + \omega^2=0$. 
The relation in (\ref{qt3}) is automatically satisfied because of the commutativity and 
cocommutativity of $H(2)$. Thus 
the $R$-matrices for $H(2)$ are in bijective correspondence with the solutions 
of the equation $2 - 2\omega + \omega^2=0$, from where we deduce that 
$R_{\pm}=1 - (1\pm i)p_{-}\ot p_{-}$ are the only quasi-triangular structures on $H(2)$.    
\end{proof}

\begin{remark}
It is not difficult to show that $H(2)_{+}=(H(2), R_{+})$ and $H(2)_{-}=(H(2), R_{-})$ 
are non-isomorphic $QT$ quasi-Hopf algebras, this means that there is no quasi-Hopf 
algebra isomorphism $\nu :\ H(2)\ra H(2)$ satisfying $(\nu \ot \nu)(R_{+})=R_{-}$. 
Indeed, if such a $\nu$ exists then 
$(1 + i)\nu(p_{-})\ot \nu (p_{-})=(1 - i)p_{-}\ot p_{-}$. If we write 
$\nu(p_{-})=ap_{-} + bp_{+}$, for some scalars $a, b\in k$, then from the above relation 
we obtain that $a^2=-i$ and $b=0$. Since $p_{\pm}^2=p_{\pm}$ and $\nu$ is an algebra 
map we get that $ap_{-}=\nu(p_{-})=\nu(p_{-}^2)=(ap_{-})^2=-ip_{-}$, and 
we conclude that $a=-i$. But $a^2=-i$, so $i\in \{-1, 0\}$, a contradiction. 
\end{remark}

To any quasi-triangular quasi-Hopf algebra $(H,R)$, we can associate a new
quasi-Hopf algebra $bos(H_0)$, called the bosonisation of $H_0$ (see \cite[Corollary 5.3]{bn2}). 
$H_0$ equals $H$ as a vector space, with a newly defined multiplication $\circ$ given
by the formula
\begin{equation}\label{ma}
h\circ h'=X^1hS(x^1X^2)\a x^2X^3_1h'S(x^3X^3_2),
\end{equation}
and left $H$-module structure given by 
$h\tr h'=h_1h'S(h_2)$, for all $h, h'\in H$. Then $bos(H_0)$ is the 
$k$-vector space $H_0\ot H$ with the following quasi-Hopf algebra structure:
\begin{eqnarray}
&&(b\times h)(b'\times h')=
(x^1\tr b)\circ (x^2h_1\tr b')\times x^3h_2h',\label{bos1}\\
&&\Delta (b\times h)\nonumber\\
&&\hspace*{5mm}
=y^1X^1\tr b_{\un{1}}\times y^2Y^1R^2x^2X^3_1h_1\ot
y^3_1Y^2R^1x^1X^2\tr b_{\un{2}}\times y^3_2Y^3x^3X^3_2h_2,\label{bos2}\\
&&\Phi _{bos(H_0)}=\b \times X^1\ot \b\times X^2\ot
\b\times X^3,\label{bos3}\\
&&s(b\times h)=(\b\times S(X^1x^1_1R^2h)\a )
(X^2x^1_2R^1\tr S_{H_0}(b)\times X^3x^2\b S(x^3)),\label{bos4}
\end{eqnarray}
for all $b, b', h, h'\in H$, where we write 
$b \times h$ and $b'\times h'$ in place of $b \ot h$  
and respectively $b'\ot h'$ to distinguish the 
new structure on $H_0\ot H$, and where 
\begin{eqnarray}
&&\un{\Delta}(b)=b_{\un{1}}\ot b_{\un{2}}:=
x^1X^1b_1g^1S(x^2R^2y^3X^3_2)
\ot x^3R^1\tr y^1X^2b_2g^2S(y^2X^3_1),\label{und}\\
&&S_{H_0}(b)=X^1R^2p^2S(q^1(X^2R^1p^1\tr b)S(q^2)X^3),\label{unant}
\end{eqnarray}
for all $b\in H$. Here $R=R^1\ot R^2$ and
$f^{-1}=g^1\ot g^2$, $p_R=p^1\ot p^2$ and $q_R=q^1\ot q^2$ are the elements 
defined by (\ref{g}) and (\ref{qr}), respectively.

The unit for $bos(H_0)$ is $\b\times 1$, the counit is 
$\va (b\times h)=\va (b)\va (h)$, for all $b, h\in H$, and the distinguished 
elements $\a$ and $\b$ are given by $\b\times \a $ and $\b\times \b$, respectively. 

Our next goal is to compute the quasi-Hopf algebra structure on $bos(H_0)$,
in the case where $H=H(2)_{+}$ or $H=H(2)_{-}$.

\begin{proposition}\prlabel{2.8}
$bos(H(2)_{+})=bos(H(2)_{-})=k[C_2\times C_2]$ as bialgebras, viewed as a 
quasi-Hopf algebra via the non-trivial reassociator 
$\Phi _x:=1 -2p^x_{-}\ot p^x_{-}\ot p^x_{-}$, where $x$ is one of the generators 
of $C_2\times C_2$, and where $p^x_{-}:=\frac{1}{2}(1 - x)$. The antipode is 
the identity map and the distinguished elements $\a$ and $\b$ are given by 
$\a=x$ and $\b=1$, respectively. In particular, $bos(H(2)_{+})=bos(H(2)_{-})$ is 
an involutory quasi-Hopf algebra. 
\end{proposition} 

\begin{proof}
Since $H(2)_{\pm}$ are commutative algebras and $\b=1$, it follows from
 (\ref{q5}) and (\ref{q6}) 
that the multiplication $\circ$ defined in (\ref{ma}) coincides with the original 
multiplication of $H(2)$. Also, from the definition of $H(2)$ it follows that 
the action $\tr$ is trivial, i.e. $h\tr h'=\va(h)h'$, for all $h, h'\in H(2)$. 
Using (\ref{bos1}) we obtain that the multiplication on $bos(H_0)$ is the 
componentwise multiplication, and by (\ref{und}) we get that the comultiplication 
$\un{\Delta}$ on $H$ reduces to 
\[
\un{\Delta}(b)=X^1b_1g^1S(y^3X^3_2)\ot y^1X^2b_2g^2S(y^2X^3_1)
=\Delta(b)(X^1X^3_2y^3\ot X^2X^3_1y^1y^2)f^{-1}.
\] 
We have that
\begin{eqnarray*}
&&(id\ot id\ot \Delta)(\Phi)=1 - 2p_{-}\ot p_{-}\ot p_{+}\ot p_{-} - 
2p_{-}\ot p_{-}\ot p_{-}\ot p_{+}, 
\end{eqnarray*}
and therefore $X^1X^3_2\ot X^2X^3_1=1$. Also, by (\ref{phih2}) we have 
\[
y^3\ot y^1y^2=1 - 2p_{-}\ot p_{-}=1\ot p_{+} + g\ot p_{-},
\]  
and a straightforward computation ensures us that the Drinfeld twist $f$ and 
its inverse $f^{-1}$ for $H(2)$ are given by 
\[
f=f^{-1}=g\ot p_{-} + 1\ot p_{+}.
\] 
Combining all these facts we get $\un{\Delta}=\Delta$, and keeping in mind 
that the action $\tr$ is trivial we conclude that the comultiplication in 
(\ref{bos2}) is the componentwise comultiplication on $H(2)\ot H(2)$. Thus 
$bos(H(2)_{+})=bos(H(2)_{-})=H(2)\ot H(2)$ as bialgebras. Hence 
$bos(H(2)_{+})=bos(H(2)_{-})$ is generated as an algebra by $x=1\times g$ and 
$y=g\times 1$, with relations $x^2=y^2=1$ and $xy=yx$. The elements 
$x$ and $y$ are grouplike elements, so $bos(H(2)_{+})=bos(H(2)_{-})=k[C_2\times C_2]$ 
as bialgebras. According to (\ref{bos3}) the reassociator of $bos(H(2)_{+})=bos(H(2)_{-})$ 
is given by  
\begin{eqnarray*} 
\Phi _x&=&1\times X^1\ot 1\times X^2\ot 1\times X^3\\
&=&1\times 1\ot 1\times 1\ot 1\times 1 - 2\times p_{-}\ot 
1\times p_{-}\ot 1\times p_{-}\\
&=& 1 - 2p_{-}^x\ot p_{-}^x\ot p_{-}^x
\end{eqnarray*} 
since $1\times p_{-}=\frac{1}{2}(1\times 1 - 1\times g)=\frac{1}{2}(1 - x)=p_{-}^x$.
Finally, 
using that $\tr$ is trivial, $\b=1$ and the axiom (\ref{q6}), we obtain $S_{H_0}=S$, 
the antipode of $H(2)$. From (\ref{bos4}) and (\ref{q6}) we deduce that 
the antipode of $bos(H(2)_{+})=bos(H(2)_{-})$ is the identity map. Clearly, 
the distinguished elements $\a$ and $\b$ are respectively
$1\times g=x$ and $1\times 1=1$. 
\end{proof}

\begin{example}\exlabel{2.9}
If $H$ and $K$ are two quasi-Hopf algebras then $H\ot K$ is also a quasi-Hopf algebra 
with the componentwise structure. In particular, the reassociator of $H\ot K$ is
\[
\Phi_{H\ot K}=(X^1_H\ot X^1_K)\ot (X^2_H\ot X^2_K)\ot (X^3_H\ot X^3_K), 
\]
where $\Phi_H=X^1_H\ot X^2_H\ot X^3_H$ and $\Phi_K=X^1_K\ot X^2_K\ot X^3_K$ 
are the reassociators of respectively $H$ and $K$. Therefore $H(2)\ot H(2)$ 
has a second quasi-Hopf algebra structure. $H(2)\ot H(2)=k[C_2\times C_2]$ as a bialgebra but 
now it is viewed as a quasi-Hopf algebra via the reassociator
\begin{eqnarray*}
&&\hspace*{-7mm}
\Phi _{x, y}=1 - 2(1\ot p_{-})\ot (1\ot p_{-})\ot (1\ot p_{-})\\
&&- 2(p_{-}\ot 1)\ot (p_{-}\ot 1)\ot (p_{-}\ot 1) + 4(p_{-}\ot p_{-})
\ot (p_{-}\ot p_{-})\ot (p_{-}\ot p_{-}).
\end{eqnarray*}
If $x=1\ot g$ and $y=g\ot 1$ are the algebra generators of $H(2)\ot H(2)$ 
then 
\[
1\ot p_{-}=\frac{1}{2}(1 - 1\ot g)=\frac{1}{2}(1 - x):=p_{-}^x,~~ 
p_{-}\ot 1=\frac{1}{2}(1 - g\ot 1)=\frac{1}{2}(1 - y):=p_{-}^y,  
\]
and $p_{-}\ot p_{-}=\frac{1}{4}(1 - g)\ot (1 - g)=\frac{1}{4}(1 - x - y + xy)
=p^x_{-}p^y_{-}$. Therefore, $\Phi _{x, y}$ comes out as 
\[
\Phi _{x, y}=(1 - 2p_{-}^x\ot p_{-}^x\ot p_{-}^x)(1 - 2p_{-}^y\ot p_{-}^y\ot p_{-}^y).
\]
Note that the distinguished elements are $\a =xy$ and $\b=1$, and that the 
antipode is the identity map. Consequently, $H(2)\ot H(2)$ is an 
involutory quasi-Hopf algebra. 
\end{example}

Another example of involutory quasi-Hopf algebra is the quantum double of $H(2)$. 

\begin{proposition}\prlabel{2.10}
The quantum double of $H(2)$ is the unital algebra generated by $X$ and $Y$ 
with relations
\[
X^2=1,~~Y^2=X,~~XY=YX.
\]
The coalgebra structure on $D(H(2))$ is given by the formulas:
\begin{eqnarray*}
&&\Delta (X)=X\ot X,~~\va(X)=1,\\
&&\Delta (Y)=-\frac{1}{2}(Y\ot Y + XY\ot Y + Y\ot XY - XY\ot XY),~~\va(Y)=-1.
\end{eqnarray*}
The reassociator, the distinguished elements $\a$ and $\b$, and the antipode are 
respectively given by
\[
\Phi _X=1 -2p_{-}^X\ot p_{-}^X\ot p_{-}^X,~~\a =X,~~\b =1,~~S(X)=X,~~S(Y)=Y,
\]
where we denoted $p_{-}^X=\frac{1}{2}(1 - X)$. Moreover, $D(H(2))$ is an involutory 
quasi-Hopf algebra.
\end{proposition}

\begin{proof}
Using the commutativity and cocommutativity of $H(2)$, (\ref{q5}), and the fact that
$\b=1$, we find that the multiplication rule (\ref{mdd}) takes the following form
on $D(H(2))$: 
\[
(\v \bowtie h)(\v'\bowtie h')=(\O ^1\O ^5\rh \v)(\O ^2\O ^4\rh \v')\bowtie \O ^3hh',
\] 
for all $\v , \v'\in H(2)^*$ and $h, h'\in H(2)$. Now, from the definition (\ref{O}) 
of $\O$ we find out that 
\[
\O ^1\O ^5\ot \O ^2\O ^4 \ot \O ^3=X^1_{(1, 1)}X^3y^1x^1f^2\ot 
X^1_{(1, 2)}X^2y^2x^2_1x^3f^1\ot X^1_2y^3x^2_2.
\]
Using the expressions of $\Phi$ and $\Phi^{-1}$ in (\ref{phih2}) we easily 
compute that 
\begin{eqnarray*}
&&X^1_{(1, 1)}X^3\ot X^1_{(1, 2)}X^2\ot X^1_2=1 -2p_{-}\ot p_{-}\ot p_{-},\\
&&x^1\ot x^2_1x^3\ot x^2_2=1 - 2p_{-}\ot p_{-}\ot p_{+},\\
&&\Phi ^{-1}(f^2\ot f^1\ot 1)=p_{-}\ot 1\ot p_{-} + p_{-}\ot g\ot p_{+} + 
p_{+}\ot 1\ot 1,
\end{eqnarray*}
where $f=g\ot p_{-} + 1\ot p_{+}$ is the Drinfeld twist of $H(2)$. By 
the above relations the multiplication of $D(H(2))$ comes out explicitly as 
\[
(\v \bowtie h)(\v '\bowtie h')=\v \v'\bowtie hh' -2(p_{-}\rh \v)(p_{-}\rh \v')\bowtie 
p_{-}hh'.
\]
Now, let $\{P_1, P_g\}$ be the dual basis of $H(2)^*$ corresponding to the basis 
$\{1, g\}$ of $H(2)$. Then $\va =P_1 + P_g$ and 
$\{\va, \mu=P_1 - P_g\}$ is clearly a basis for $H(2)^*$. Now let $X=\va \bowtie g$ and 
$Y=\mu \bowtie 1$. Since $p_{-}\rh P_1=\frac{1}{2}\mu$ and $p_{-}\rh P_g=
-\frac{1}{2}\mu$, we obtain that
\begin{eqnarray*}
&&X^2=(\va \bowtie g)(\va \bowtie g)=\va \bowtie 1 - 2(p_{-}\rh \va)(p_{-}\rh \va)=1,\\
&&XY=YX=\mu\bowtie g,\\
&&Y^2=\mu^2\bowtie 1 - 2(p_{-}\rh \mu)^2\bowtie p_{-}
=\va \bowtie 1 - \va \bowtie 2p_{-}=\va \bowtie g=X, 
\end{eqnarray*}
which are the multiplication rules that we stated. A computation as in the 
proof of \prref{2.8} shows that the reassociator of $D(H(2))$ has the desired form.

Since $H(2)$ is commutative and cocommutative, $\b=1$ and 
$\Phi ^{-1}=\Phi =Y^2\ot Y^1\ot Y^3$, by (\ref{cddf}) we have 
\[
\Delta(\v \bowtie h)=p^1_1p^2\rh \v_2\bowtie p^1_2X^1h_1\ot X^2_1X^3\rh 
\v_1\bowtie X^2_2h_2,
\]
for all $\v\in H(2)^*$ and $h\in H(2)$. On the other hand,
\begin{eqnarray*}
&&\hspace*{-2cm}
p^1_1p^2\ot p^1_2X^1\ot X^2_1X^3\ot X^2_2\\
&=&(p_{+}\ot 1\ot 1\ot 1 - p_{-}\ot g\ot 1\ot 1)(1 - 2\ot p_{-}\ot p_{-}\ot p_{+})\\
&=&p_{+}\ot 1\ot 1\ot 1 - p_{-}\ot g\ot 1\ot 1 - 2\ot p_{-}\ot p_{-}\ot p_{+}.
\end{eqnarray*}
We then have 
\[
\Delta (X)=\Delta(\va\bowtie g)=\va \bowtie g\ot \va\bowtie g=X\ot X,  
\]
and since $\Delta (P_1)=P_1\ot P_1 + P_g\ot P_g$ and $\Delta(P_g)=P_1\ot P_g + P_g\ot P_1$ 
we get that $\Delta(\mu)=\Delta (P_1 - P_g)=(P_1 - P_g)\ot (P_1 - P_g)=\mu \ot \mu$, and 
therefore 
\begin{eqnarray*}&&\hspace*{-10mm}
\Delta(Y)=p_{+}\rh \mu \bowtie 1\ot \mu \bowtie 1
- p_{-}\rh \mu \bowtie g \ot \mu \bowtie 1
-2\mu \bowtie p_{-}\ot p_{-}\rh \mu \bowtie p_{+}\\
&=&- XY\ot Y  - \frac{1}{2}(Y - XY)\ot (Y + XY)\\
&=&-\frac{1}{2}(Y\ot Y + XY\ot Y + Y\ot XY - XY\ot XY),
\end{eqnarray*}
as needed. It follows from (\ref{coddf}) that $\va(X)=1$ and $\va (Y)=-1$.

Finally, in our particular situation (\ref{anddf}) takes the form  
\[
S(\v \bowtie h)=p^1_1p^2U^1f^2\rh \v \bowtie p^1_2U^2f^1h
\]
for all $\v\in H(2)^*$ and $h\in H$. But $p^1_1p^2f^2\ot p^1_2f^1=U^1\ot U^2=g\ot 1$, 
so the antipode for $D(H(2))$ is the identity map. Obviously, $\a =\va \bowtie g=X$, 
$\b =\va \bowtie 1=1$, and our proof is complete. 
\end{proof}

\begin{remark}
Since $p^1_1p^2\ot p^1_2=p_{+}\ot 1 - p_{-}\ot g$ we have from (\ref{rdd}) that the canonical 
$R$-matrix for $D(H(2))$ is $R=p_{+}^X\ot 1 - p_{-}^X\ot XY$, where, as usual, 
$p_{\pm}^X=\frac{1}{2}(1\pm X)$. 
\end{remark}

At first sight there is no relationship between $D(H(2))$ and $H(2)\ot H(2)$,
so it comes as a surprise that these two Hopf algebras are twisted equivalent.
To show this, we will use the structure of the quantum double associated to a factorizable 
quasi-Hopf algebra, see \cite{bt0}.

Recall from \cite{bt0} that a $QT$ quasi-Hopf algebra $(H, R)$ is called factorizable if 
the $k$-linear map ${\cal Q}:\ H^*\ra H$ given for all $\chi \in H^*$ by
\begin{equation}\label{qf1}
{\cal Q}(\chi)=\big\le \chi , S(X^2_2\tilde{p}^2)f^1R^2r^1U^1X^3\big\ri 
X^1S(X^2_1\tilde{p}^1)f^2R^1r^2U^2,
\end{equation}
is bijective. Here $r^1\ot r^2$ is another copy of $R$, and 
$U=U^1\ot U^2$, $f^{-1}=g^1\ot g^2$ and $p_L=\tilde{p}^1\ot \tilde{p}^2$ 
are the elements defined by (\ref{U}), (\ref{g}) and (\ref{ql}), respectively.  
The first step is to prove that $H(2)$ is a factorizable quasi-Hopf algebra.

\begin{proposition}
For $(H(2), R)$ with $R$ as in \prref{2.6} the map ${\cal Q}$ from (\ref{qf1}) 
has the following form for all $\chi\in H(2)^*$
\[
{\cal Q}(\chi)=\chi(1)p_{-} + \chi(g)p_{+}.
\]
Since $\{p_{-}, p_{+}\}$ and $\{1, g\}$ are bases for $H(2)$ it 
follows that ${\cal Q}$ is bijective, so $H(2)$ is factorizable.
\end{proposition}

\begin{proof}
For $H(2)$ the element $p_L$ has the form 
\[
p_L=X^1X^2\ot X^3=1 - 2p_{-}\ot p_{-}=1 - (1 - g)\ot p_{-}=1\ot p_{+} + g\ot p_{-}=f.
\]
Also, we can easily see that $X^1X^2_1\ot X^2_2X^3=1$ and since $f$ is an involution 
we conclude that 
\[
X^2_2X^3\tilde{p}^2f^1\ot X^1X^2_1\tilde{p}^1f^2=1.
\]
On the other hand, since $\omega ^2 -2\omega =-2$ it follows that 
$R^2r^1\ot R^1r^2=(1 -\omega p_{-}\ot p_{-})^2=1 - 2p_{-}\ot p_{-}$. We have already 
seen that $U=g\ot 1$, and therefore 
\begin{eqnarray*}
&&\hspace*{-2cm}
S(X^2_2\tilde{p}^2)f^1R^2r^1U^1X^3\ot X^1S(X^2_1\tilde{p}^1)f^2R^1r^2U^2\\
&&=(1 - 2p_{-}\ot p_{-})(g\ot 1)=1\ot p_{-} + g\ot p_{+}.
\end{eqnarray*}
It is now clear that 
${\cal Q}(\chi)=\chi(1)p_{-} + \chi(g)p_{+}$, for all $\chi\in H(2)^*$, and this 
finishes the proof.  
\end{proof}  

The structure of the quantum double of a finite dimensional factorizable quasi-Hopf algebra 
$(H, R)$ can be found in \cite[Theorem 5.4]{bt0}. More precisely, since $(H, R)$ is 
quasi-triangular there exist two quasi-Hopf algebra morphisms 
$\pi , \tilde{\pi}:\ D(H)\ra H$ covering the natural inclusion $i_D:\ H\ra D(H)$. They are given 
by the formulas
\[
\pi(\v \bowtie h)=\v (q^2R^1)q^1R^2h~~{\rm and}~~
\tilde{\pi}(\v \bowtie h)=\v (q^2\ov{R}^2)q^1\ov{R}^1h, 
\]
for all $\v \in H^*$ and $h\in H$, where $R^{-1}=\ov{R}^1\ot \ov{R}^2$. 
Now, if we define 
\begin{equation}\label{psptwist}
\mathbf{F}=Y^1_1x^1X^1y^1_1\ot Y^1_2x^2R^2X^3y^2
\ot Y^2x^3R^1X^2y^1_2\ot Y^3y^3,
\end{equation}
and $\mathbf{U}=\ov{R}^1g^2\ot \ov{R}^2g^1$, then $F$ is a twist on $H\ot H$, $\mathbf{U}$ 
is an invertible element of $H\ot H$ 
and $\zeta :\ D(H)\ra (H\ot H)_{\mathbf{F}}^{\mathbf{U}}$, 
given  by 
$\zeta ({\bf D})=\widetilde{\pi }({\bf D}_1)\ot \pi ({\bf D}_2)$, 
for all ${\bf D}\in D(H)$,
is a quasi-Hopf algebra isomorphism. We make this result explicit for the
quasitriangular quasi-Hopf algebra $(H(2), R)$.

\begin{proposition}\prlabel{2.13}
Let $X, Y$ be the algebra generators of $D(H(2))$ defined in \prref{2.10}, and $x, y$ 
the generators of $H(2)\ot H(2)\cong k[C_2\times C_2]$, the quasi-Hopf algebra 
described in \exref{2.9}. Let $\omega _{\pm}=1\pm i$  and consider the elements 
\begin{eqnarray*}
&&\mathbf{U}_{\pm}=p_{+}\ot 1 + p_{-}\ot g + \omega _{\pm}p_{-}\ot p_{-},\\
&&\mathbf{F}_{\pm}=1 - 2p_{-}^xp_{-}^y\ot p_{-}^xp_{+}^y 
- 2p_{+}^xp_{-}^y\ot p_{-}^xp_{-}^y -
\omega _{\pm}p_{-}^x\ot p_{-}^y, 
\end{eqnarray*}  
where $p_{\pm}^x=\frac{1}{2}(1\pm x)$ and $p_{\pm}^y=\frac{1}{2}(1\pm y)$. 
Then the maps $\zeta _{\pm}:\ D(H(2))\ra (H(2)\ot H(2))_{\mathbf{F}_{\pm}}^{\mathbf{U}_{\pm}}$, 
given by 
\[
\zeta_{\pm}(X)=xy,~~\zeta_{\pm}(Y)=-1 + \omega_{\pm}p_{-}^x + \omega_{\mp}p_{-}^y
=-{1\over 2}(\omega_{\pm} x+\omega_{\mp}y),
\]
are quasi-Hopf algebra isomorphisms. 
\end{proposition}

\begin{proof}
Since $H(2)$ is a factorizable quasi-Hopf algebra 
everything will follow from the general isomorphism presented above. For $H(2)$ 
we have $q_R=1\ot p_{+} - g\ot p_{-}$. Also, it is easy to see that the inverse 
of $R_{\pm}=1 -\omega_{\pm}p_{-}\ot p_{-}$ is 
$R_{\mp}=1 -\omega_{\mp}p_{-}\ot p_{-}$, and therefore
\begin{eqnarray*}
&&q^2R^1_{\pm}\ot q^1R^2_{\pm}=p_{+}\ot 1 - p_{-}\ot g - \omega_{\pm}p_{-}\ot p_{-},\\
&&q^2\ov{R}^2_{\pm}\ot q^1\ov{R}^1_{\pm}=q^2R^1_{\mp}\ot q^1R^2_{\mp}=
p_{+}\ot 1 - p_{-}\ot g - \omega_{\mp}p_{-}\ot p_{-}.
\end{eqnarray*}   
From the structure of $D(H(2))$ in \prref{2.10} we compute that
$\pi(X)=\tilde{\pi}(X)=g$,
\[
\pi(Y)=\pi(\mu \bowtie 1)=\mu(p_{+})1 - \mu(p_{-})g - \omega_{\pm}\mu(p_{-})p_{-}
=-g - \omega_{\pm}p_{-},
\]
and, in a similar way, $\tilde{\pi}(Y)=-g - \omega_{\mp}p_{-}$. We get that 
$\pi(XY)=-1 + \omega_{\pm}p_{-}$ and $\tilde{\pi}(XY)=-1 + \omega_{\mp}p_{-}$, 
so $\zeta _{\pm}(X)=g\ot g=xy$ and 
\[
\zeta_{\pm}(Y)=-\frac{1}{2}(\pi\ot \tilde{\pi})\left(Y\ot Y + XY\ot Y 
+ Y\ot XY - XY\ot XY\right).
\]
After some straightforward computations we obtain 
\begin{eqnarray*}
&&\pi(Y)\ot \tilde{\pi}(Y)=xy + 2p_{-}^xp_{-}^y + \omega_{\pm}xp_{-}^y + 
\omega_{\mp}yp_{-}^x,\\
&&\pi(XY)\ot \tilde{\pi}(Y)=x - 2p_{-}^xp_{-}^y -\omega_{\pm}xp_{-}^y +\omega_{\mp}p_{-}^x,\\
&&\pi(Y)\ot \tilde{\pi}(XY)=y - 2p_{-}^xp_{-}^y +\omega_{\pm}p_{-}^y -\omega_{\mp}yp_{-}^x,\\
&&\pi(XY)\ot \tilde{\pi}(XY)=1 + 2p_{-}^xp_{-}^y - \omega_{\pm}p_{-}^y -\omega_{\mp}p_{-}^x.
\end{eqnarray*}
Thus, we can compute:
\begin{eqnarray*}
\zeta_{\pm}(Y)&=&-\frac{1}{2}(xy + x + y - 1 - 4p_{-}^xp_{-}^y + 2\omega_{\pm}p_{-}^y + 
2\omega_{\mp}p_{-}^x)\\
&=&1 - x - y -\omega_{\pm}p_{-}^y - \omega_{\mp}p_{-}^x\\
&=&-1 + (2 -\omega_{\mp})p_{-}^x + (2 - \omega{\pm})p_{-}^y
=-1 + \omega_{\pm}p_{-}^x + \omega_{\mp}p_{-}^y,
\end{eqnarray*}
and this is exactly what we need. Finally, one can easily see that 
the corresponding elements $\mathbf{U}_{\pm}$ and $\mathbf{F}_{\pm}$ for 
$(H(2), R_{\pm})$ are exactly the ones defined in the statement, 
we leave the details to the reader. 
\end{proof}

\begin{remarks}\relabel{2.14}
1) Keeping the notation used in \prref{2.13}, we have that 
\[
(\zeta_{\pm}\ot\zeta_{\pm}\ot\zeta_{\pm})(\Phi _X)
=\Phi_{xy}:=1 - 2p_{-}^{xy}\ot p_{-}^{xy}\ot p_{-}^{xy}, 
\]
where $\Phi_X$ is the reassociator of $D(H(2))$ and $p_{-}^{xy}:=\frac{1}{2}(1 - xy)$. 
We then have that $\Phi_{xy}$ is a $3$-cocycle for $k[C_2\times C_2]$ and 
$\Phi_{xy}=(\Phi_{x,y})_{\mathbf{F}}=(\Phi_x\Phi_y)_{\mathbf{F}}$, because 
of (\ref{g2}). Here $\Phi_y=1 - 2p_{-}^y\ot p_{-}^y\ot p_{-}^y$ is the $3$ cocycle 
on $k[C_2\times C_2]$ corresponding to $y$. In other words we have proved that 
the $3$-cocycles $\Phi_{xy}$ and $\Phi_{x}\Phi_{y}$ are equivalent. 

2) It follows from \prref{2.13} that $k[C_4]$ and $k[C_2\times C_2]$ are isomorphic as
algebras if char$(k)\neq 2$ and $k$ contains a primitive fourth root of 1. This is well-known
and can be easily seen directly: $k[C_4]$ and $k[C_2\times C_2]$ are isomorphic
(as Hopf algebras even) to their duals and the two duals are both isomorphic to $k^4$
as algebras. More explicitly, we have the following: $\zeta_+=\gamma\circ\beta\circ\alpha$,
where $\alpha,\beta,\gamma$ are the following three algebra isomorphisms.
$\{e_1,e_2,e_3,e_4\}$ is the standard basis of $k^4$.
\begin{eqnarray*}
&&
\alpha:\ k[C_4]=k[Y]/(Y^4-1)\to k^4,~~\alpha(Y)=e_1+ie_2-e_3-ie_4;\\
&&
\beta:\ k^4\to k[C_2\times C_2]=k[x,y],\\
&&\hspace*{2cm}
\beta (e_1)=p_{+}^x p_{+}^y,~~\beta (e_2)=p_{-}^x p_{+}^y,~~
\beta (e_3)=p_{+}^x p_{-}^y;~~
\beta (e_4)=p_{-}^x p_{-}^y;\\
&&
\gamma:\ k[C_2\times C_2]\to k[C_2\times C_2],~~\gamma(x)=xy,~~
\gamma(y)=-x,~~\gamma(xy)=-y.
\end{eqnarray*}
\end{remarks}

\section{The pivotal structure of 
${}_H{\cal M}^{\rm fd}$ when $H$ is involutory}\selabel{pivotal}
\setcounter{equation}{0}
If ${\cal C}$ is a monoidal category with left duality, then the functor $(-)^{**}:\
{\cal C}\to {\cal C}$ is strongly monoidal: we have an isomorphism $\phi_0:\
I\to I^{**}$, and
for $V,W\in {\cal C}$, we have the
following family of isomorphisms in ${\cal C}$:
\[
\phi_{V, W}:\ V^{**}\ot W^{**}~~\stackrel{\l_{W^*, V^*}}{\longrightarrow}~~
(W^*\ot V^*)^*~~\stackrel{(\l^{-1}_{V, W})^*}{\longrightarrow}~~(V\ot W)^{**}.
\]
$\l_{V, W}:\ W^*\ot V^*\ra (V\ot W)^*$ is the isomorphism described in 
\cite[Proposition XIV.2.2]{k}), and $(\l^{-1}_{V, W})^*$ is the transpose in ${\cal C}$ of the 
morphism $\l^{-1}_{V, W}$, see \cite[XIV.2]{k}.
 
By definition, a pivotal structure on ${\cal C}$ is an isomorphism $i$ between the strongly 
monoidal functors ${\rm Id}$ and $(-)^{**}$. This means that $i$ is a natural 
transformation satisfying the coherence conditions
\begin{equation}\eqlabel{coherence}
\phi_{V, W}\circ (i_V\ot i_W)=i_{V\ot W},~~\phi_0=i_I.
\end{equation}

The importance of pivotal structures lies in the following fundamental result of
Etingof, Nikshych and Ostrik \cite[Prop. 8.24 and 8.23]{eno}. 
Recall first that the category  ${\cal C}={}_H{\cal M}^{\rm fd}$ of finite dimensional 
left modules over a quasi-Hopf algebra $H$ is monoidal with left duality. The structure 
is the following. If $U, V, W$ are left $H$-modules, define 
$a_{U, V, W}: (U\ot V)\ot W\ra U\ot (V\ot W)$ by 
\[
a_{U, V, W}((u\ot v)\ot w)=\Phi\cdot (u\ot (v\ot w)).
\]
Then ${}_H{\cal M}$ becomes a monoidal category with tensor product $\ot$ 
given via $\Delta$, associativity constraints $a_{U, V, W}$, unit $k$ as a trivial 
$H$-module and the usual left and right unit constraints. In addition,  
every object $V$ of ${\cal C}$ has a (left) dual object $V^*$, the 
linear dual of $V$, with left $H$-action $\le h\cdot \v , v\ri =\le \v ,S(h)\cdot v\ri$. 
The evaluation and coevaluation maps are given  by the formulas
\[
{\rm ev}_V(\v \ot v)=\v (\a \cdot v),~~
{\rm coev}_V(1)=\sum\limits_i\b \cdot v_i\ot v^i,
\]
for all $\v\in V^*$ and $v\in V$.
Here $\{v_i\}_i$ is a basis of $V$ with dual basis $\{v^i\}_i$ in $V^*$.

The next result was proved in \cite[Prop. 8.24 and 2.3]{eno}. 

\begin{theorem}\thlabel{piv}
If $H$ is a semisimple quasi-Hopf algebra over an algebraically closed field of characteristic zero 
then ${}_H{\cal M}^{\rm fd}$ has a unique pivotal structure 
$j:\ {\rm Id}\to (-)^{**}$ such that for any simple object $V$ of ${\cal C}$, 
$\un{\rm dim}_k(V)={\rm dim}_k(V)$, where 
$\un{\rm dim}_k(V):={\rm ev}_{V^*}\circ (j_V\ot {\rm id}_{V^*})\circ {\rm coev}_V$
is the (categorical) dimension of $V$ in ${\cal C}$.
\end{theorem}  

The aim of this Section is to compute explicitly this unique pivotal structure 
for an involutory quasi-Hopf algebra $H$. Note that, due to the trace formula 
proved in \cite{bt}, $H$ is semisimple, so that \thref{piv} can be applied.

We first give a description of all pivotal structures on ${}_H{\cal M}^{\rm fd}$, when $H$ 
is finite dimensional. 
For this, we do not need any assumption on the groundfield $k$.

\begin{proposition}\prlabel{3.1}
Let $H$ be a finite dimensional quasi-Hopf algebra over a field $k$. Then we 
have a bijective correspondence between pivotal structures $i$ on ${\cal C}
={}_H{\cal M}^{\rm fd}$ 
and invertible elements ${\mf g}_i\in H$ satisfying 
\begin{equation}\label{prop1}
S^2(h)={\mf g}_i^{-1}h{\mf g}_i,
\end{equation}
for all $h\in H$ and 
\begin{equation}\label{pivstr2}
\Delta({\mf g}_i)=({\mf g}_i\ot {\mf g}_i)(S\ot S)(f_{21}^{-1})f, 
\end{equation} 
where $f=f^1\ot f^2$ is the Drinfeld twist defined in (\ref{f}) and 
$f_{21}=f^2\ot f^1$.  
\end{proposition}

\begin{proof}
Since $H$ is finite dimensional its antipode $S$ is bijective, cf. \cite{bc1}. 

Let $V\in {}_H{\cal M}^{\rm fd}$. We have a $k$-linear isomorphism
$V\to V^{**}$, and $V^{**}$ can be regarded as $V$ with newly defined left
$H$-action $h\cdot v=S^2(h)v$.

It has been pointed out in the literature (see \cite{mn,sch}), that there is a
bijective correspondence between natural isomorphisms between the
functors Id and $(-)^{**}$ and invertible elements ${\mf g}_i\in H$ satisfying
(\ref{prop1}). Let us sketch this correspondence. Let $i:\ {\rm Id}\to (-)^{**}$
be a natural isomorphism, and let ${\mf g}_i= S^2(i_H^{-1}(1))$, ${\mf k}_i=
i_H(1)$. For all $h\in H$, we have $i_H(h)=S^2(h){\mf k}_i$ and
$i_H^{-1}(h)=S^{-2}(h{\mf g}_i)$. In particular, $1=i_H^{-1}({\mf k}_i)=
S^{-2}({\mf k}_i{\mf g}_i)$ and $1=i_H(S^{-2}({\mf g}_i))={\mf g}_i{\mf k}_i$,
hence ${\mf k}_i={\mf g}_i^{-1}$. Now take $V\in {}_H{\cal M}^{\rm fd}$, and
fix $v\in V$. From the naturality of $i$, we deduce that
$i_V(v)=(i_V\circ f)(1)=(f\circ i_H)(1)= {\mf g}_i^{-1}v$. This means that $i$
is completely determined by ${\mf g}_i$:
\begin{equation}\label{pivstr1}
i_V(v)= {\mf g}_i^{-1} v.
\end{equation}
Now take $V=H$ and $v=h$. (\ref{pivstr1}) tells us that $S^2(h){\mf g}_i^{-1}=
i_H(h)={\mf g}_i^{-1}h$, and it follows that (\ref{prop1}) is satisfied.

For $H$ a quasi-Hopf algebra and ${\cal C}={}_H{\cal M}^{\rm fd}$ the isomorphisms 
$\l_{V, W}$ were computed in \cite[Proposition 4.2]{bcp}, namely
\[
\lambda_{V, W}(w^*\ot v^*)(v\ot w)=\le v^*,f^1\cd v\ri \le w^*,f^2\cd w\ri,    
\]
for all $v\in V$, $w\in W$, $v^*\in V^*$ and $w^*\in W^*$, where $f=f^1\ot f^2$ 
is the Drinfeld twist defined in (\ref{f}). Observe that the isomorphism
$\lambda_{V,W}$ is denoted $\phi^*_{W,V}$ in \cite{bcp}.
It is not difficult to see at this point 
that $\l^{-1}_{V, W}$ is given by 
\[
\l^{-1}_{V, W}(\psi )=\sum\limits_{i, j}\psi (g^1\cdot v_i\ot g^2\cdot w_j)w^j\ot v^i,
\]
for all $\psi\in (V\ot W)^*$, where $\{v_i\}_i$ and $\{v^i\}_i$ are dual bases of $V$ and $V^*$, 
$\{w_j\}_j$ and $\{w^j\}_j$ are dual bases of $W$ and $W^*$, and $f^{-1}=g^1\ot g^2$ is the 
inverse of the Drinfeld twist.

It is then easy to establish that the first condition from \equref{coherence} is equivalent 
to the following equivalent conditions:
\begin{eqnarray*}
&&\hspace*{-1cm} \phi_{V, W}(i_V(v)\ot i_W(w))=i_{V\ot W}(v\ot w)\\
&\Leftrightarrow& 
\l_{W^*, V^*}(i_V(v)\ot i_W(w))\circ \l_{V, W}^{-1}=i_{V\ot W}(v\ot w)\\
&\Leftrightarrow& 
\l_{W^*, V^*}(i_V(v)\ot i_W(w))(w^*\ot v^*)=i_{V\ot W}(v\ot w)(\l_{V, W}(w^*\ot v^*))\\
&\Leftrightarrow&
i_W(w)(f^1\cdot w^*)i_V(v)(f^2\cdot v^*)=\l_{V, W}(w^*\ot v^*)({\mf g}_i^{-1}\cdot (v\ot w))\\
&\Leftrightarrow&
v^*(S(f^2){\mf g}_i^{-1}\cdot v)w^*(S(f^1){\mf g}_i^{-1}\cdot w)=
v^*(f^1({\mf g}_i^{-1})_1\cdot v)w^*(f^2({\mf g}_i^{-1})_2\cdot w),
\end{eqnarray*} 
for all $V, W\in {\cal C}$, $v\in V$, $w\in W$, $v^*\in V^*$ and $w^*\in W^*$. Since $H$ 
is an object of ${\cal C}$, this last condition is equivalent to 
\[
\Delta({\mf g}_i^{-1})=f^{-1}(S\ot S)(f_{21})({\mf g}_i^{-1}\ot {\mf g}_i^{-1}),
\]
which is equivalent to (\ref{pivstr2}). Now take $x\in k$. It follows from (\ref{pivstr1})
that $j_k(x)=\varepsilon(g_i^{-1})x$. Since $\phi_0:\ k\to k^{**}$ is the identity, we see
that the second condition from \equref{coherence} is equivalent to $\varepsilon({\mf g}_i)=1$.
If (\ref{pivstr2}) is satisfied, then $\varepsilon({\mf g}_i)^2=\va ({\mf g}_i)$ (apply $\varepsilon\ot
\varepsilon$ to (\ref{pivstr2})), and it follows that $\varepsilon({\mf g}_i)=1$. This completes
our proof.
\end{proof}

The categorical dimension corresponding to $i$ is now given by the formula
\begin{equation}\label{prop2}
\un{\rm dim}_k(V)=\sum\limits_iv^i({\mf g}_i^{-1}\b S(\a)\cdot v_i)
=\sum\limits_iv^i({\mf g}_iS(\b)\a\cdot v_i),
\end{equation}
for all $V\in {\cal C}$.

The element ${\mf g}$ corresponding to the unique pivotal structure 
$j$ in \thref{piv} was computed in \cite{mn}. Since $H$ is semisimple there exists a 
(left and right) integral $\Lambda$ in $H$ such that $\va(\Lambda)=1$, cf. \cite{p}. 
By \cite[Corollary 8.5]{mn} we then have
\[
{\mf g}=q^2\Lambda_2p^2S(q^1\Lambda_1p^1),
\]
where $p_R=p^1\ot p^2$ and $q_R=q^1\ot q^2$ are the elements defined in (\ref{qr}). 
We should note that the above formula for ${\mf g}$ can be immediately 
obtained from the trace formula proved in \cite{bt}. An alternative way to 
compute ${\mf g}$ can be found in \cite[Lemma 3.1]{sch}, in the case where $\b$
is invertible.

If $H$ is an involutory quasi-Hopf algebra, then $S^2$ is inner, and induced by
$\b S(\a)$. We will now show that the invertible element ${\mf g}$ corresponding to
unique pivotal structure on ${\cal C}$ is precisely
$\b S(\a)$. In particular, ${\mf g}^{-1}=S(\b)\a$, and the equalities 
in (\ref{prop2}) become trivial. 

\begin{proposition}\prlabel{pivforinv}
Let $H$ be a finite dimensional involutory quasi-Hopf algebra over an algebraic closed field of 
characteristic zero. Then the element ${\mf g}$ corresponding to the pivotal structure $j$ 
in \thref{piv} is equal to $\b S(\a)$. 
\end{proposition}

\begin{proof}
By \cite[Lemma 2.1]{bc1}, we have for any left integral $t$ in $H$ that 
\[
t_1\ot t_2=\b q^1t_1\ot q^2t_2=q^1t_1\ot S^{-1}(\b)q^2t_2.
\]
Replacing $H$ by $H^{\rm op}$ we find for any right integral $r$ in $H$ that, 
\[
r_1\ot r_2=r_1p^1S^{-1}(\a)\ot r_2p^2=r_1p^1\ot r_2p^2\a .
\] 
Now, since $\Lambda$ is both a left and right integral in $H$ we compute:
\begin{eqnarray*}
&&\hspace*{-2cm}
{\mf g}=q^2\Lambda _2p^2S(q^1\Lambda ^1p^1)=S^{-1}(\b^{-1})\Lambda_2p^2S(\Lambda_1p^1)
=S^{-1}(\b^{-1})\Lambda_2\a ^{-1}S(\Lambda_1)\\
\hspace*{1cm}&{{\rm (\ref{sp})}\atop =}&
\va(\Lambda)S^{-1}(\b^{-1})\a^{-1}=\b S(\a)S(\b^{-1})S(\b)\a \a^{-1}
=\b S(\a).
\end{eqnarray*}
In the penultimate equality, we used  (\ref{inv}) in its equivalent form, 
$S^{-1}(h)=\b S(\a)$ $S(h)S(\b)\a$, for all $h\in H$. 
\end{proof}

Assume that $H$ is an involutory quasi-Hopf algebra. It is a natural question to ask
whether the Drinfeld double $D(H)$ is also involutory. We will present a sufficient
condition in \prref{2.15}.
In order to simplify the computations 
we need the following formulas
\begin{eqnarray}
&&f^1_1p^1\ot f^1_2p^2S(f^2)=g^1S(\tilde{q}^2)\ot g^2S(\tilde{q}^1),\label{for1}\\
&&S(U^1)\tilde{q}^1U^2_1\ot \tilde{q}^2U^2_2=f,\label{for2}
\end{eqnarray}
where $p_R=p^1\ot p^2$ and 
$q_L=\tilde{q}^1\ot \tilde{q}^2$ are the elements defined 
in (\ref{qr}) and (\ref{ql}), and $f=f^1\ot f^2$ is the Drinfeld's twist defined 
in (\ref{f}) with its inverse $f^{-1}=g^1\ot g^2$ as in (\ref{g}). (\ref{for1},\ref{for2})
follow easily from the axioms and the basic properties of a quasi-Hopf algebra.

\begin{proposition}\prlabel{2.15}
Let $H$ be an involutory quasi-Hopf algebra such that 
\begin{equation}\label{for3}
\Delta(S(\b)\a)=f^{-1}(S\ot S)(f_{21})(S(\b)\a \ot S(\b)\a),
\end{equation}
where $f_{21}=f^2\ot f^1$. Then $D(H)$ is an involutory quasi-Hopf algebra. 
\end{proposition}

\begin{proof}
For all $\v\in H^*$ and $h\in H$ we compute:
\begin{eqnarray*}
&&\hspace*{-2.5cm}
S^2_D(\v \bowtie h)\\
&=&S_D(p^1_1U^1\rh \ov{S}^{-1}(\v)\lh f^2\smi (p^2)
\bowtie p^1_2U^2)(\va \bowtie S(S(h)f^1))\\
&{{\rm (\ref{anddf})}\atop =}&
(\va\bowtie  S(p^1_2U^2)F^1)(P^1_1{\cal U}^1\rh 
\ov{S}^{-1}(p^1_1U^1\rh \ov{S}^{-1}(\v)\lh f^2\smi(p^2))\\
&&\hspace*{-15mm}
\lh F^2\smi(P^2)\bowtie P^1_2{\cal U}^2)
(\va \bowtie S(f^1))(\va \bowtie S^2(h))\\
&{{\rm (\ref{mdd}, \ref{for1})}\atop =}&
(S(p^1_2U^2)_1g^1S(\tilde{q}^2))_1{\cal U}^1p^2S(f^2)\rh \ov{S}^{-2}(\v)\\
&&\hspace*{-15mm}
\lh S(p^1_1U^1)\tilde{q}^1\smi(S(p^1_2U^2)_2g^2)\bowtie 
(S(p^1_2U^2)_1g^1S(\tilde{q}^2))_2{\cal U}^2S(f^1))
(\va \bowtie S^2(h))\\
&{{\rm (\ref{ca}, \ref{ql1a})}\atop =}&
(g^1_1S(p^1\tilde{q}^2U^2_2)_1{\cal U}^1p^2S(f^2)\rh \ov{S}^{-2}(\v)\lh 
S(U^1)\tilde{q}^1U^2_1\smi(g^2)\\
&&\hspace*{-15mm}
\bowtie g^1_2S(p^1\tilde{q}^2U^2_2)_2{\cal U}^2S(f^1))(\va \bowtie S^2(h))\\
&{{\rm (\ref{U}, \ref{ca}, \ref{qr1})}\atop =}&
(g^1_1G^1S(f^2\tilde{q}^2_2U^2_{(2, 2)})\rh \ov{S}^{-2}(\v)\lh 
S(U^1)\tilde{q}^1U^2_1\smi (g^2)\\
&&\hspace*{-15mm}
\bowtie g^1_2G^2S(f^1\tilde{q}^2_1U^2_{(2, 1)}))(\va \bowtie S^2(h))\\
&{{\rm (\ref{for2})}\atop =}&
g^1_1G^1S(f^2F^2_2)\rh \ov{S}^{-2}(\v)\lh F^1\smi(g^2)\bowtie g^1_2G^2S(f^1F^2_1)S^2(h),
\end{eqnarray*}
where ${\cal U}^1\ot {\cal U}^2$, $F^1\ot F^2$ and $G^1\ot G^2$ are second 
copies of $U$, $f$ and $f^{-1}$, respectively.

Using (\ref{for3}) twice and (\ref{ca}) we obtain that
\begin{eqnarray*}
&&(\Delta\ot id)(\Delta(S(\b)\a))\\
&&\hspace*{1cm}
=g^1_1G^1S(f^2F^2_2)S(\b)\a \ot 
g^1_2G^2S(f^1F^2_1)S(\b)\a \ot g^2S(F^1)S(\b)\a .  
\end{eqnarray*} 
By (\ref{inv}) and (\ref{inva}) we have 
\[
S(\b S(\a))=S^2(\a)S(\b)=S(\b)\a ^2\b S(\b \a)=S(\b)\a ^2\a ^{-1}=S(\b)\a,
\] 
or, equivalently, $\smi (S(\b)\a)=\b S(\a)$. We then have 
\begin{eqnarray*}
&&\hspace*{-1cm}
(\va \bowtie S(\b)\a)(\v \bowtie h)(\va \bowtie \a S(\b))\\
&{{\rm (\ref{mdd})}\atop =}&
g^1_1G^1S(f^2F^2_2)S(\b)\a\rh \v \lh \b S(\a)\smi(g^2S(F^1))\bowtie 
g^1_2G^2S(F^1f^2_1)S^2(h)\\
&=&g^1_1G^1S(f^2F^2_2)\rh \ov{S}^{-2}(\v)\lh F^1\smi(g^2)
\bowtie g^1_2G^2S(f^1F^2_1)S^2(h)\\
&=&S^2_D(\v \bowtie h),
\end{eqnarray*}
for all $\v \in H^*$ and $h\in H$, and this means that $D(H)$ is an involutory 
quasi-Hopf algebra. 
\end{proof}

\begin{remarks}
1) The formula (\ref{for3}) holds for any finite dimensional involutory 
quasi-Hopf algebra over an algebraic closed field of characteristic zero. Indeed, 
by \prref{pivforinv} we have ${\mf g}^{-1}=S(\b)\a$, so (\ref{for3}) follows from 
(\ref{pivstr2}). Moreover, we believe that (\ref{for3}) is satisfied for an arbitrary involutory 
quasi-Hopf algebra; this would imply that $H$ is involutory if and only if $D(H)$ is 
involutory. Somehow this should follow naturally from the equality 
$\Delta(S^2(h))=\Delta(S(\b)\a h\b S(\a))$ which is equivalent to 
\begin{eqnarray*}
&&\hspace*{-5mm}
\Delta(S(\b)\a)\Delta(h)\Delta(\b S(\a))\\
&&\hspace*{5mm}
=f^{-1}(S\ot S)(f_{21})(S(\b)\a \ot S(\b)\a)\Delta(h)(\b S(\a)\ot \b S(\a))
(S\ot S)(f^{-1}_{21})f.
\end{eqnarray*}

2) For $H(2)$, the condition (\ref{for3}) reduces to $\Delta(g)=g\ot g$ which is just 
part of the definition of $H(2)$. So \prref{2.15} gives a new proof for the fact that
$D(H(2))$ is involutory.

3) One of the aims in \cite{mn} was to compute the Frobenius-Schur indicator 
for the irreducible representations of $D^{\omega}(G)$. Note that $D^{\omega}(G)$ 
is a particular case of the quasi-Hopf algebra $D^{\omega}(H)$ roughly described 
in \exsref{exsinv} 1). Thus $D^{\omega}(G)$ is an involutory quasi-Hopf algebra, so 
the element ${\mf g}$ corresponding to $D^{\omega}(G)$ is $\beta$, and its inverse is 
$\b^{-1}=S(\b)$.   
\end{remarks}

\section{Representations of involutory quasi-Hopf algebras}\selabel{3}
\setcounter{equation}{0}
The goal of this Section is to study the representations of an involutory quasi-Hopf 
algebra $H$ over a field $k$. In the case where $H$ is semisimple, we will show that the 
characteristic of $k$ does not divide the dimension of any finite dimensional 
absolutely simple $H$-module. 
Following \cite[Definition 3.42]{cr} a left $H$-module $V$ is called absolutely simple if for 
every field extension $k\subseteq K$, $K\ot V$ is a simple $K\ot H$-module or, equivalently, 
if every $H$-endomorphism of $V$ is on the form $c~id_V$ for some scalar $c\in k$ (see 
\cite[Theorem 3.43]{cr}).

The case when $H$ is not semisimple is treated as well. In this case the characteristic 
of $k$ divides the dimension of any finite dimensional projective $H$-module.

In the case of Hopf algebras, the first result was initially proved for Hopf algebras by Larson in \cite{l} 
and afterwards by Lorenz \cite{lz} in a slightly more general shape.
The Hopf algebra version of the second result is also due to Lorenz \cite{lz}.

Before we are able to prove these results over involutory quasi-Hopf algebras,
we need some preliminary results.

Let $V$ and $W$ be left $H$-modules and denote by ${\rm Hom}_H(V, W)$ the set of $H$-linear 
morphisms between $V$ and $W$. If $V$ is finite dimensional then  
${\rm Hom}_H(V, W)\cong {\rm Hom}_H(k, W\ot V^*)$ 
as $k$-vector spaces. Actually, this isomorphism works in any monoidal 
category ${\cal C}$, in the sense that for any object $V$ of ${\cal C}$ admitting a left dual object $V^*$ and 
for any $W\in {\cal C}$ we have ${\rm Hom}_{\cal C}(V, W)\cong {\rm Hom}_{\cal C}(\un{1}, W\ot V^*)$, where $\un{1}$ is 
the unit object of ${\cal C}$. Indeed, it is well-known that 
\[
{\rm Hom}_{\cal C}(V, W)\ni \vartheta \mapsto \left(\un{1}\rTo^{{\rm coev}_V}V\ot V^*\rTo^{\vartheta \ot {\rm Id}_{V^*}}
W\ot V^*\right)\in {\rm Hom}_{\cal C}(\un{1}, W\ot V^*)
\]
is an isomorphism in ${\cal C}$; its inverse is given by   
\begin{eqnarray*}
&&\hspace*{-1cm}
{\rm Hom}_{\cal C}(\un{1}, W\ot V^*)\ni \varsigma \mapsto 
\left(V\rTo^{l^{-1}_V}\un{1}\ot V\rTo^{\varsigma\ot {\rm Id}_V}(W\ot V^*)\ot V\right. \\
&&\hspace{5mm}
\left.
\rTo^{a_{W, V^*, V}}W\ot (V^*\ot V)
\rTo^{{\rm Id}_W\ot {\rm ev}_V}W\ot \un{1}\rTo^{r_W}W\right)\in {\rm Hom}_{\cal C}(V, W),
\end{eqnarray*}
where $a, l, r$ are the associativity, and the left and right unit constraints of ${\cal C}$. 

In particular, if $V\in {}_H{\cal M}^{fd}$ then ${\rm End}_H(V):={\rm Hom}_H(V, V)$ 
is isomorphic to ${\rm Hom}_H(k, V\ot V^*)$, as $k$-vector space. 

Secondly, for any morphism $a: V\ra V^{**}$ in ${}_H{\cal M}^{\rm fd}$ the categorical trace 
$\un{\rm Tr}_V(a)$ of $a$ is defined as the scalar corresponding to 
\[
k\rTo^{{\rm coev}_V}V\ot V^*\rTo^{a\ot{\rm Id}_{V^*}}V^{**}\ot V^*\rTo^{{\rm ev}_{V^*}}k.
\]

We are now able to prove one of the main results of this Section.

\begin{theorem}\label{te3.4}
Let $H$ be a semisimple involutory quasi-Hopf algebra over a field $k$ of 
characteristic $p\geq 0$. Then $p$ does not divide the dimension of any 
finite dimensional absolutely simple $H$-module.
\end{theorem}

\begin{proof}
Since $H$ is involutory we have $S^2(h)={\mf g}^{-1}h{\mf g}$, for all $h\in H$, where 
${\mf g}=\b S(\a)$ and ${\mf g}^{-1}=S(\b)\a$ is its inverse in $H$. Then for any 
$V\in {}_H{\cal M}^{\rm fd}$ we get that $a: V\ra V^{**}$ defined by $a(v)(v^*)=v^*({\mf g}^{-1}\cdot v)$, for 
all $v\in V$ and $v^*\in V^*$, is an isomorphism in ${}_H{\cal M}^{\rm fd}$. Moreover, 
its categorical trace $\un{\rm Tr}_V(a)$ coincides with the usual $k$-dimension of $V$. Indeed, 
if $\{v_i\}_{i=\ov{1, n}}$ is a basis in $V$ with dual basis 
$\{v^i\}_{i=\ov{1, n}}$ in $V^*$ we then have 
\[
\un{\rm Tr}_V(a)=\sum\limits_{i=1}^n a(\b\cdot v_i)(\a\cdot v^i)=\sum\limits_{i=1}^n (\a\cdot v^i)({\mf g}^{-1}\b\cdot v_i)
=\sum\limits_{i=1}^n v^i(S(\b\a)\a\b\cdot v_i)={\rm dim}_k(V),
\]  
because of \leref{2.2}. 

Assume now that $V$ is a finite dimensional absolutely simple $H$-module. Thus $V\in {}_H{\cal M}^{fd}$ 
and ${\rm dim}_k({\rm End}_H(V))=1$. We claim that $\un{\rm Tr}_V(a)\not=0$; this 
would imply that ${\rm dim}_k(V)\not=0$, as required. 

By the way of contradiction, if $\un{\rm Tr}_V(a)=0$ we then obtain a sequence 
\[
0\ra k\rTo^{{\rm coev}_V}V\ot V^*\rTo^{{\rm ev}_{V^*}(a\ot{\rm Id}_{V^*})}k\ra 0 
\]
in ${}_H{\cal M}^{fd}$ with ${\rm Im}({\rm coev}_V)\subseteq {\rm Ker}({\rm ev}_{V^*}(a\ot{\rm Id}_{V^*}))$. 
Since $H$ is semisimple $k$ is a projective left $H$-module, and so we get a sequence in ${}_k{\cal M}$,  
\[
0\ra k \rTo^{\zeta}{} {\rm Hom}_H(k, V\ot V^*)\rTo^{\upsilon} k \ra 0, 
\]
with $\zeta$ and $\upsilon$ $k$-linear morphisms satisfying 
${\rm Im}(\zeta)\subseteq {\rm Ker}(\upsilon)$. Now, if $\psi: k\ra {\rm Hom}_H(k, V\ot V^*)$ 
is a $k$-linear map such that $\upsilon\psi ={\rm Id}_k$ and $\theta:=\psi\upsilon$ then 
$\theta^2=\theta$, hence ${\rm Hom}_H(k, V\ot V^*)={\rm Im}(\theta)\oplus {\rm Ker}(\theta)$. 
Taking into account that ${\rm Im}(\theta)={\rm Im}(\psi)\cong k$ and that 
${\rm Ker}(\theta)={\rm Ker}(\upsilon)\supseteq {\rm Im}(\zeta)\cong k$ it follows that 
${\rm dim}_k({\rm Hom}_H(k, V\ot V^*))\geq 2$. Thus 
$1={\rm dim}_k({\rm End}_H(V))={\rm dim}_k({\rm Hom}_H(k, V\ot V^*))\geq 2$,
a contradiction.  
\end{proof} 

\begin{remark}
Some arguments in the preceding proof have also been used in \cite{eno} in the context 
of fusion categories. 
\end{remark}

Using Schur's Lemma and Theorem \ref{te3.4}, we immediately obtain the following
result.

\begin{corollary}
Let $H$ be a semisimple involutory quasi-Hopf algebra over an algebraically closed 
field of characteristic $p\geq 0$. Then $p$ does not divide the dimension of any 
finite dimensional simple $H$-module.
\end{corollary}

We will now focus attention to the case where $H$ is not semisimple. 
 Let $V$ and $W$ be two left $H$-modules. Then ${\rm Hom}_k(V, W)$ is a left
$H$-module, with structure given by the formula
$(h\cd \psi )(v)=h_1\cd \psi (S(h_2)\cd v)$,
for all $h\in H$, $\psi\in {\rm Hom}_k(V, W)$ and $v\in V$. Furthermore, 
if $V$ is finite dimensional then ${\rm Hom}_k(V, W)\cong W\ot V^*$ as left 
$H$-modules. To see this take $\{v_i\}_{i=\ov{1, n}}$ a basis in $V$ with dual basis 
$\{v^i\}_{i=\ov{1, n}}$ and define  
\[
\xi: {\rm Hom}_k(V, W)\ni \chi\mapsto \sum\limits_{i=1}^n\chi(v_i)\ot v^i\in W\ot V^*.
\]
Then $\xi$ is a left $H$-linear isomorphism with 
\[
\xi^{-1}: W\ot V^*\ni w\ot v^*\mapsto \left(v\mapsto v^*(v)w\right)\in {\rm Hom}_k(V, W).
\]
(The verification of all these details is left to the reader.) In particular, 
if $V$ is a finite dimensional left $H$-module then 
${\rm End}_k(V):={\rm Hom}_k(V, V)\cong V\ot V^*$, as left $H$-modules. 

In order to simplify 
the proof of Theorem \ref{te3.6}, we first need the following result.

\begin{proposition}\label{pr3.5}
Let $H$ be a finite dimensional quasi-Hopf algebra and $P, Q$ finite dimensional 
projective left $H$-modules. Then:
\begin{itemize}
\item[(i)] $P^*={\rm Hom}_k(P, k)$ is a projective left $H$-module;
\item[(ii)] $P\ot Q$ is a projective left $H$-module, where the $H$-module 
structure of $P\ot Q$ is defined by the comultiplication $\Delta$ of $H$;
\item[(iii)] ${\rm End}_k(P)$ is a projective left $H$-module.
\end{itemize}  
\end{proposition}

\begin{proof}
(i) Note that, if $V$ is a left $H$-module then $V^*$, the linear dual space of V, 
is a left $H$-module via the $H$-action $(h\cd v^*)(v)=v^*(S(h)\cd v)$, for all 
$v^*\in V^*$, $h\in H$ and $v\in V$.\\
${\;\;\;}$
Since $P$ is finite dimensional it follows that $P$ is a finitely generated projective left 
$H$-module. Therefore, there exists  $n\in \NN$ and a left $H$-module $P'$ such that 
$P\oplus P'\cong H^n$ as left $H$-modules. Thus
\begin{eqnarray*}
&&\hspace*{-1cm}
(H^n)^*={\rm Hom}_k(H^n ,k)\cong {\rm Hom}_k(P\oplus P', k)\\
&&\hspace*{2cm}
\cong {\rm Hom}_k(P, k)\oplus 
{\rm Hom}_k(P', k)=P^*\oplus P'^*,
\end{eqnarray*}
as left $H$-modules. Now, $H$ is finite dimensional, so from the proof of 
\cite[Theorem 4.3]{hn3} we know that the map 
\[
H\ni h\mapsto \left(h'\mapsto \l (h'S(h))\right)\in H^*
\]
is bijective. Here $\l$ is a non-zero left cointegral on $H$, we refer to 
\cite{hn3} for the definition of a left cointegral. 
Replacing $H$ by $H^{\rm op, cop}$ we get that $H\cong H^*$ as left $H$-modules.

Then we obtain that 
$H^n\cong (H^*)^n\cong (H^n)^*\cong P^*\oplus P'^*$, 
as left $H$-modules, so $P^*$ is a projective left $H$-module.

(ii) We follow the same line as above. There exist $n,m\in \NN$
and two left $H$-modules $P'$ and $Q'$ such that 
$P\oplus P'\cong H^n$ and $Q\ot Q'\cong H^m$, as left $H$-modules. 
We then have  
\[
(H\ot H)^{nm}\cong H^n\ot H^m\cong (P\ot Q)\oplus (P\ot Q')\oplus (P'\ot Q)\oplus 
(P'\ot Q'),
\]
as left $H$-modules. It now suffices to show that $H\ot H$, with the diagonal
$H$-action, is free as a left $H$-module. To this end, 
we will show that the map   
\[
\mu : {}_{\cd}H\ot {}_{\cd}H\ra {}_{\cd}H\ot H,~~
\mu (h\ot h')=\tqlb h'_2\ot \smi (\tqla h'_1)h,
\]
is a left $H$-linear isomorphism. Here we denoted by 
${}_{\cd}H\ot {}_{\cd}H$ and ${}_{\cd}H\ot H$ the $k$-vector space $H\ot H$,
respectively with the diagonal left $H$-action, and left $H$-action given by
left multiplication. $q_L=\tqla \ot \tqlb$ 
is the element defined in (\ref{ql}).
For all $h, h', h''\in H$ we have that
\begin{eqnarray*}
&&\hspace*{-2cm}
\mu (h''\cd (h\ot h'))=\mu (h''_1h\ot h''_2h')
=\tqlb h''_{(2, 2)}h'_2\ot \smi (\tqla h''_{(2, 1)}h'_1)h''_1h\\
&\equal{(\ref{ql1a})}&h''\tqlb h'_2\ot \smi (\tqla h'_1)h
=h''\mu(h\ot h'),
\end{eqnarray*}
proving that $\mu$ is left $H$-linear. It is easy check that 
the map 
\[
\mu ^{-1}:\ {}_{\cd}H\ot H\ra {}_{\cd}H\ot {}_{\cd}H,~~
\mu ^{-1}(h\ot h')=h_1\tpla h'\ot h_2\tplb 
\]
is the inverse of $\mu$. More precisely, (\ref{ql1}) and (\ref{pqla}) imply that 
$\mu ^{-1}\circ \mu =id$, while (\ref{ql1a}) and (\ref{pql}) imply that 
$\mu \circ \mu ^{-1}=id$, we leave the verification of the details to the reader.

(iii)
If $P$ is finite dimensional, then ${\rm End}_k(P)\cong P\ot P^*$ 
as left $H$-modules. (iii) is now an immediate application of 
(i) and (ii).
\end{proof}

We are now able to prove the second important result of this Section. 

\begin{theorem}\label{te3.6}
Let $H$ be a finite dimensional involutory quasi-Hopf algebra 
over a field of characteristic $p\geq 0$. If $H$ is not semisimple, 
then $p$ divides the dimension of any 
finite dimensional projective $H$-module. 
\end{theorem}

\begin{proof}
Let $P$ be a finite dimensional 
projective left $H$-module and suppose that $p$ does not divide 
${\rm dim}_k(P)$. If $a: P\ra P^{**}$ is the isomorphism in ${}_H{\cal M}^{fd}$ 
defined in the proof of Theorem \ref{te3.4}, specialized for $V=P$, then 
$\un{\rm Tr}_P(a)={\rm dim}_k(P)\not=0$ in $k$, and so 
${\rm coev}_P: k\ra P\ot P^*$ is a split monomorphism.  
Hence, $k$ is isomorphic to a direct summand of $P\ot P^*$ which is a projective 
left $H$-module by Proposition \ref{pr3.5}. It follows that $k$ is a projective left $H$-module.

Now, since $\va$ is left $H$-linear and surjective, we obtain that there 
exists an $H$-linear map $\vartheta :\ k\ra H$ such that $\va \circ \vartheta =id_k$. Then 
$t=\vartheta (1_k)$ is a left integral in $H$ since 
\[
ht=h\vartheta (1_k)=\vartheta (h\cd 1_k)=\va (h)\vartheta (1_k)=\va (h)t,
\]
for all $h\in H$. Moreover, $\va (t)=\va (\vartheta (1_k))=1_k$, 
so by the Maschke-type Theorem proved in \cite{p} we obtain that 
$H$ is semisimple, a contradiction.         
\end{proof}
\section{Involutory dual quasi-Hopf algebras with non-zero integrals}\selabel{4}
\setcounter{equation}{0}
Let $H$ be a finite dimensional quasi-Hopf algebra over a field of characteristic 
$p\geq 0$. It was proved in \cite{bt} that 
\[{\rm Tr}\left(h\mapsto \b S(\a)S^2(h)S(\b)\a \right)=\va(r)\l (\smi (\a)\b),
\]
where $r$ is a non-zero right integral in $H$ (this means a left integral in $H^{\rm op}$) 
and $\l$ is a non-zero left cointegral on $H$ such that $\l (S(r))=1$. 
In particular, if $H$ is involutory it follows from \leref{2.2} that $S(\a)S^2(h)S(\b)\a=h$, and 
\[
{\rm dim}_k(H)={\rm Tr}(id_H)=\va (r)\l (\smi (\a)\b).
\]
Thus a finite dimensional involutory quasi-Hopf algebra is both 
semisimple and cosemisimple if and only if ${\rm dim}_k(H)\not=0$ in $k$. 
We recall from \cite{hn3, bt} 
that $H$ is called cosemisimple if there is a left cointegral on $H$ 
such that $\l (\smi (\a)\b)\not=0$. Consequently, we obtain that a finite dimensional 
involutory quasi-Hopf algebra over a field of characteristic zero is always semisimple 
and cosemisimple.

By duality, we obtain that a finite dimensional involutory dual quasi-Hopf algebra over a 
field of characteristic zero is both semisimple and cosemisimple. The aim of this Section is 
to study the infinite dimenensional case. In fact, we will prove that an involutory 
co-Frobenius dual quasi-Hopf algebra over a field of characteristic zero is cosemisimple. 
Our approach is based on the methods developed in \cite{dnt}.

Throughout this Section, $A$ will be a dual quasi-Hopf algebra. 
Following \cite{maj}, a dual quasi-bialgebra $A$ is a
coassociative coalgebra $A$ with comultiplication
$\Delta $ and counit $\va $ together with coalgebra morphisms
$m_A :\ A\ot A\ra A$ (the multiplication;
we write $m_A(a\ot b)=ab$)
and $\eta _A:\ k\ra A$ (the unit; we write $\eta _A(1)=1$), and an
invertible element $\v \in (A\ot A\ot A)^*$ (the reassociator),
such that for all $a, b, c, d\in A$ the following relations
hold (summation understood):
\begin{eqnarray}
&&a_1(b_1c_1)\v (a_2, b_2, c_2)=
\v (a_1, b_1, c_1)(a_2b_2)c_2,\label{dq1}\\
&&
1a=a1=a,\label{dq2}\\
&&\v (a_1, b_1, c_1d_1)\v (a_2b_2, c_2,
d_2)=\v (b_1, c_1, d_1)\v (a_1, b_2c_2, d_2)
\v (a_2, b_3, c_3),\label{dq3}\\
&&\v (a, 1, b)=\va (a)\va (b)\label{dq4}.
\end{eqnarray}

$A$ is called a dual quasi-Hopf algebra if, moreover,
there exist an anti-morphism $S$ of the coalgebra $A$ and elements
$\a , \b \in H^*$ such that, for all $a\in A$:
\begin{eqnarray}
&&S(a_1)\a (a_2)a_3=\a (a)1,
\mbox{${\;\;\;}$}
a_1\b (a_2)S(a_3)=\b (a)1,\label{dq5}\\
&&\v (a_1\b (a_2),
S(a_3), \a (a_4)a_5)=
\v ^{-1}(S(a_1), \a (a_2)a_3, \b (a_4)S(a_5))=\va (a).\label{dq6}
\end{eqnarray}

It follows from the axioms that $S(1)=1$ and $\a (1)\b (1)=1$,
so we can assume that $\a (1)=\b (1)=1$. Moreover (\ref{dq3}) and
(\ref{dq4}) imply
\begin{equation}\label{dq7}
\v (1, a, b)=\v (a, b, 1)=\va (a)\va (b),~~\forall ~~a, b\in A.
\end{equation}

Note that, if $A$ is a dual quasi-bialgebra then $A^*$, the linear dual space of $A$, 
is an algebra with multiplication given by the convolution and unit $\va$.

We call a dual quasi-Hopf algebra $A$ 
involutory if 
\begin{equation}\label{dinv}
S^2(a)=\b (S(a_1))\a (a_2)a_3\b (a_4)\a (S(a_5)),~~\forall ~~a\in A.
\end{equation}
The proof of the following result is formally dual to the proof of 
\leref{2.2} and \prref{2.4}, and is left to the reader.

\begin{proposition}\label{pr4.1}
Let $A$ be an involutory dual quasi-Hopf algebra. Then $(\b \circ S)\a$ is 
convolution invertible with $((\b \circ S)\a )^{-1}=\b (\a \circ S)$. In particular, 
the square of the antipode is coinner, so the antipode $S$ is bijective.
 
Moreover, the elements $\a$ and $\b$ are convolution invertible and for all 
$a\in A$ the following relations hold: 
\[
S(a _3)\a ^{-1}(a_2)a_1=\a ^{-1}(a)1~~{\rm and}~~
a_3\b ^{-1}(a_2)S(a_1)=\b ^{-1}(a)1.
\]
\end{proposition}

Let $A$ be a dual quasi-Hopf algebra. The connection between integrals and 
the ideal $A^{*{\rm rat}}$ was given in \cite{bc1}. $A^{*{\rm rat}}$ is our notation for the 
left rational part of $A^*$, and coincides with the right 
rational part of $A^*$, see \cite{bc1}. 
Also recall that a left integral on $A$ is an element $T\in A^*$ 
such that $a^*T=a^*(1)T$, for all $a^*\in A^*$, and that $\int _l$ is the standard 
notation for the space of left integrals on $A$. Finally, note that $A^{*{\rm rat}}\neq 0$ 
if and only if $\int _l\not=0$, if and only if $A$ is a left or right co-Frobenius coalgebra. 
In this case ${\rm dim}_k(\int _l)=1$ (see \cite{bc1}).

Now let $\sigma :\ A\ot A\ra A^*$ be defined by $\sigma (a\ot b)(c)=\v (c, a, b)$, 
for all $a, b, c\in A$. $\sigma$ is convolution invertible, with 
inverse given by $\sigma ^{-1}(a\ot b)(c)=\v ^{-1}(c, a, b)$, for all $a, b, c\in A$. 
Now define $\theta ^*:\  \int _l\ot A\ra A^{*{\rm rat}}$ by
\[
\theta ^*(T\ot a)=\sigma (S(a_5)\ot \a(a_6)a_7)(T\leftharpoondown a_4)
\sigma ^{-1}(S(a_3)\ot \b(S(a_2))S^2(a_1)), 
\]
for all $T\in \int_l$ and $a\in A$.
The right $A$-action $\leftharpoondown$ on $A^*$ is given by 
$(a^*\leftharpoondown a)(b)=\le a^*, bS(a)\ri$, for all $a^*\in A^*$ and $a, b\in A$. 
It is proved in \cite[Proposition 4.2]{bc1} that $\theta^*$ is a well-defined  
isomorphism of right $A$-comodules.

Assume that $S$ is bijective, and define the elements $p_R, q_R\in (A\ot A)^*$ by
\begin{equation}\label{dpqr}
p_R(a, b)=\v ^{-1}(a, b_1, S(b_3))\b (b_2),~~
q_R(a, b)=\v (a, b_3, \smi (b_1))\a (\smi (b_2)).
\end{equation}  
Then the map $\theta^*$ can be written in the following form:
\[
\theta ^*(T\ot a)(b)=q_R(b_1, S(a_3))T(b_2S(a_2))p_R(b_3, S(a_1)),
\]
for all $a, b\in A$ and $T\in \int _l$.

Let $A$ be a co-Frobenius dual quasi-Hopf algebra with non-zero left integral $T$.
For every $a^*\in A^{*{\rm rat}}$, there exists $a\in A$ such that
\begin{equation}\label{imp}
a^*(b)=\omega _T(b, S(a)),
\end{equation}
for all $b\in A$. The map $\omega _T:\ A\ot A\ra k$ is defined by the fomrula
\begin{equation}\label{om}
\omega _T(b, a)=q_R(b_1, a_1)T(b_2a_2)p_R(b_3, a_3),
\end{equation}
for all $a,b\in A$.
The next two Lemmas will be crucial in the sequel.

\begin{lemma}\label{le4.2}
Let $A$ be an involutory dual quasi-Hopf algebra, and $T$ a left integral on $A$.
The map  
$\omega _T$ satisfies the following formula, for all $a\in A$: 
\begin{equation}\label{prom}
\omega _T(a_2, S(a_1))=T(1)\b (S(a_1))\a (a_2).
\end{equation} 
\end{lemma}

\begin{proof}
For all $a,b\in A$, we have
\begin{equation}\label{simp}
q_R(a_1, S(b_2))\b (b_3)T(a_2S(b_1))=T(aS(b)).
\end{equation}
This formula is  the formal dual of the first equality in 
\cite[Lemma 2.1, (2.3)]{bc1}.\\
Also, from (\ref{dinv}) and Proposition \ref{pr4.1} we have for all $a\in A$ that 
\begin{equation}\label{sdinv}
a=\b (a_1)\a (S(a_2))S^2(a_3)\b (S(a_4))\a (a_5).
\end{equation}
We then compute, for all $a\in A$ that
\begin{eqnarray*}
&&\hspace*{-2cm}
\omega _T(a_2, S(a_1))\equal{(\ref{sdinv})}
\b (a_2)\a (S(a_3))\omega _T(S^2(a_4), S(a_1))\b (S(a_5))\a (a_6)\\
&\equal{(\ref{om})}&q_R(S^2(a_6), S(a_3))\b (a_4)\a (S(a_5))
T(S^2(a_7)S(a_2))\\
&&\hspace*{1cm}\times 
p_R(S^2(a_8), S(a_1))\b (S(a_9))\a (a_{10})\\
&\equal{(\ref{simp})}&
\a (S(a_3))T(S^2(a_4)S(a_2))p_R(S^2(a_5), S(a_1))\b (S(a_6))\a (a_7)\\
&\equal{(\ref{dq5})}&
T(1)\a (S(a_2))p_R(S^2(a_3), S(a_1))\b (S(a_4))\a (a_5)\\
&\equal{(\ref{dpqr}, \ref{dq6})}&
T(1)\va (S(a_1))\b (S(a_2))\a (a_3)
=T(1)\b (S(a_1))\a (a_2),
\end{eqnarray*} 
as claimed, and this completes the proof. 
\end{proof}

Our next result is a generalization of \cite[Lemma 1]{dnt}.

\begin{lemma}\label{le4.3}
Let $H$ be an involutory dual quasi-Hopf algebra with a non-zero 
left integral $T$ and $J$ a finite dimensional non-zero left 
$A$-subcomodule of $A$ which is a direct summand of $A$ as a left $A$-comodule. 
Then ${\rm dim}_k(J)=cT(1)$, for some scalar $c\in k$.    
\end{lemma}

\begin{proof}
Let $J'$ be a left $A$-subcomodule of $A$ such that $J\oplus J'=A$, and define 
$a^*\in A^*$ by $a^*(a+a')=\varepsilon(a)$, for $a\in J$ and $a'\in J'$.
Since $J$ is finite dimensional it follows that ${\rm Ker}(a^*)$ contains 
a left $A$-subcomodule of $A$ of finite codimension. It then follows from \cite[Corollary 2.2.16]{dnr} 
that $a^*\in A^{*{\rm rat}}$, and we have 
an element $a\in A$ such that $a^*(b)=\omega _T(b, S(a))$, for all $b\in A$.

For all $b,c\in A$, we have the following formula (cf. \cite[(4.16)]{bc1})
\[
q_R(b_2, S(c_2))T(b_3S(c_1))b_1=q_R(b_1, S(c_2))T(b_2S(c_1))c_3,
\]
Using this formula and the definition 
of $\omega _T$ we compute that 
\begin{eqnarray*}
&&\hspace*{-15mm} a^*\rh b=a^*(b_2)b_1=\omega _T(b_2, S(a))b_1=
q_R(b_2, S(a_3))T(b_3S(a_2))p_R(b_4, S(a_1))b_1\\
&=&
q_R(b_1, S(a_3))T(b_2S(a_2))p_R(b_3, S(a_1))a_4=\omega _T(b, S(a_1))a_2,
\end{eqnarray*}
for all $b\in B$. 
Now consider a basis $\{a_i\}_{i=\ov{1, n}}$ for $J$ and $\{a'_{\l}\}_{\l \in \Lambda }$ 
a basis for $J'$ and then write 
\[
\Delta (a)=\sum \limits _{i=1}^n b_i\ot a_i + \sum \limits _{\l \in \Lambda}b'_{\l}\ot a'_{\l}
\]  
for some elements $b_i, b'_{\l}\in A$. We then have 
\begin{eqnarray*}
&&\hspace*{-2cm}
a_j=\va ((a_j)_2)(a_j)_1=a^*((a_j)_2)(a_j)_1=a^*\rh a_j
=\omega _T(a_j, S(a_1))a_2\\
&=&\sum \limits _{i=1}^n \omega _T(a_j, S(b_i))a_i + \sum \limits _{\l \in \Lambda}
\omega _T(a_j, S(b'_{\l}))a'_{\l},
\end{eqnarray*} 
for any $j\in\{1,\cdots,n\}$, so 
\[
\omega _T(a_j, S(b_i))=\d _{i, j}~~{\rm and}~~
\omega _T(a_j, S(b'_{\l}))=0,
\]
for all $i, j\in \{1, \cdots , n\}$ and $\l \in \Lambda$, where $\d _{i, j}$ is the 
Kronecker's symbol.

In a similar way, we compute for all $\l '\in \Lambda$ that 
\begin{eqnarray*}
&&\hspace*{-2cm}
0=a^*((a'_{\l '})_2)(a'_{\l '})_1=a^*\rh a'_{\l '}=\omega _T(a'_{\l '}, S(a_1))a_2\\
&=& \sum \limits _{i=1}^n\omega _T(a'_{\l '}, S(b_i))a_i + \sum \limits _{\l \in \Lambda}
\omega _T(a'_{\l '}, S(b'_{\l}))a'_{\l}
\end{eqnarray*}
We find that
$\omega _T(a'_{\l '}, S(b_i))=\omega _T(a'_{\l '}, S(b'_{\l}))=0$, for all 
$i\in \{1, \cdots , n\}$ and $\l , \l '\in \Lambda$. It follows now that 
\[
{\rm dim}_k(J)=\sum \limits _{i=1}^n\omega _T(a_i, S(b_i)) + \sum \limits _{\l \in \Lambda}
\omega _T(a'_{\l}, S(b'_{\l}))=\omega _T(a_2, S(a_1)).
\]
From Lemma \ref{le4.2} we conclude that ${\rm dim}_k(J)=T(1)\b (S(a_1))\a (a_2)$, so 
${\rm dim}_k(J)=c T(1)$ for $c=\b (S(a_1))\a (a_2)$, as claimed.  
\end{proof}

We can now prove the main result of this Section, generalizing \cite[Theorem 2]{dnt}.
The remaining arguments are a purely coalgebraic flavour, and are identical to the
arguments in \cite{dnt}.

\begin{theorem}\label{te4.4}
Let $A$ be an involutory dual quasi-Hopf algebra with non-zero integral over a field 
of characteristic zero. Then $A$ is cosemisimple. 
\end{theorem}

\begin{proof}
We have that $A$ is an injective left $A$-comodule, so there exists 
an injective envelope $J$ of $k1_A$ such that $J\subseteq A$. Being injective, 
$J$ is a direct summand of $A$. Since 
$A$ is co-Frobenius by \cite[Theorem 3]{lin} we obtain that $J$ is finite dimensional. 
Applying Lemma \ref{le4.3} we deduce that $T(1)\not=0$ and from 
\cite[Theorem 4.10]{bc1} we conclude that $A$ is cosemisimple.  
\end{proof}


\end{document}